\documentclass{amsart}

\usepackage[american, british]{babel}
\usepackage{hyperref, tikz}

\newtheorem{theorem}{Theorem}[section]
\newtheorem{proposition}[theorem]{Proposition}
\newtheorem{corollary}[theorem]{Corollary}
\newtheorem{lemma}[theorem]{Lemma}
\newtheorem{example}[theorem]{Example}
\theoremstyle{remark}
\newtheorem{remark}[theorem]{Remark}

\numberwithin{equation}{section}

\DeclareMathOperator{\dv}{div}
\DeclareMathOperator{\vl}{vol}
\DeclareMathOperator{\hs}{Hess}
\DeclareMathOperator{\e}{e}
\DeclareMathOperator{\im}{Im}
\DeclareMathOperator{\supp}{supp}

\begin{document}

\title[Isoparamtric functions and mean curvature with navigation]{Isoparametric functions and mean curvature in manifolds with Zermelo navigation}

\author[Benigno Oliveira Alves]{Benigno Oliveira Alves}
\address{Benigno Oliveira Alves \textup{(corresponding author)} \hfill\break\indent Instituto de Matem\'{a}tica e Estat\'{\i}stica, Universidade Federal da Bahia \hfill\break\indent Rua Bar\~ao de Jeremoabo, 40170-115 Salvador, Bahia, Brazil}
\email{benignoalves@ufba.br; gguialves@hotmail.com}
\thanks{Warm thanks to Professors Marcos Alexandrino, Miguel {\'A}ngel Javaloyes, and Miguel Dom{\'i}nguez V{\'a}zquez for the helpful discussions.}

\author[Patr\'icia Mar\c cal]{Patr\'icia Mar\c cal}
\address{Patr\'icia Mar\c cal}
\email{patriciaqmarcal@gmail.com}

\date{\today}


\begin{abstract}
The generalized Zermelo navigation problem looks for the shortest time paths in an environment, modeled by a Finsler manifold $(M, F)$, under the influence of wind or current, represented by a vector field $W$. The main objective of this paper is to investigate the relationship between the isoparametric functions on the manifold $M$ with and without the presence of the vector field $W$. Our work generalizes results in \cite{he2020classifications, dong2020isoparametric, xu2020some, xu2021isoparametric, he2022isoparametric}. For the positive-definite cases, we also compare the mean curvatures in the manifold. Overall, we follow a coordinate-free approach.
\end{abstract}

\maketitle
\bibliographystyle{amsplain}

\section{Introduction}

The starting point for the study of isoparametric hypersurfaces is a conceptually simple question in geometric optics: What are the waves whose speed remains constant on each wavefront? Formal investigations in $\mathbb{R}^3$ date back to the early 20th century \cite{laura1918sopra, somigliana1918sulle, segre1924proprieta}. Laura showed that the possible wavefronts are strongly restricted \cite{laura1918sopra}. Somigliana proved that the wavefronts must be either parallel planes, concentric spheres, or coaxial cylinders  \cite{somigliana1918sulle}. Segre later obtained independently the same result \cite{segre1924proprieta}. So did Levi-Civita in the late 1930s \cite{levi1937famiglie}. Soon after, Segre and Cartan considered the problem in Riemannian manifolds \cite{segre1938famiglie, cartan1938familles}.

A nonconstant real-valued function $f$ on a Riemannian manifold $M$ is called \emph{transnormal} if there is a real smooth function $b$ such that
\begin{equation}
\label{transnRiem}
\vert \nabla f \vert^2 = b \circ f .
\end{equation}
A transnormal function is \emph{isoparametric} if there is a real continuous function $a$ with
\begin{equation}
\label{isoRiem}
\triangle f = a \circ f ,
\end{equation}
where $\triangle$ denotes the Laplace operator on $M$; in this case, the regular level sets of $f$ form a family of isoparametric hypersurfaces -- the wavefronts. Condition~\ref{transnRiem} means the wavefronts are parallel to each other, which can be interpreted as Huygens' principle \cite{somigliana1918sulle}; Equation~\ref{isoRiem} then indicates each regular level set has constant mean curvature \cite{cartan1938familles}. Segre classified isoparametric hypersurfaces in Euclidean spaces \cite{segre1938famiglie}; Cartan, isoparametric hypersurfaces in hyperbolic spaces \cite{cartan1938familles}. The classification in spheres started with Cartan \cite{cartan1938familles, cartan1939sur}, but it was concluded only recently \cite{chi2020isoparametric}, with the contribution of many mathematicians. For a survey on the topic until the beginning 21st century, see \cite{thorbergsson1999survey} and the references therein.

The theory of isoparametric hypersurfaces in Finsler spaces has just begun, and there are several reasons that might explain its delay. On one hand, there is no canonical volume form associated with the Finsler structure; the most important choices are the Busemann-Hausdorff and the Holmes-Thompson volume forms. Moreover, the volume form associated with the induced Finsler structure in a hypersurface need not coincide with the volume form induced to the hypersurface by a unit normal vector field \cite{shen1997curvature}. In fact, there are several definitions for ``induced volume form''; see related discussion in \cite{chakerian1996integral, schneider1997integral}. On the other hand, there are different notions for the Laplacian \cite{anastasiei1993absolute, antonelli1993stochastic, bao1996hodge, bao1996eigenforms, bellettini1996anisotropic, shen1997curvature, shen1998nonlinear}, and both the gradient and the most natural choice for the Laplacian are nonlinear operators that cannot be explicitly expressed in the general case \cite{shen1997curvature}. On top of that, since the nonlinear Laplacian $\triangle f$ is undefined at points where $\mathrm{d} f = 0$, many methods of Riemannian geometry cannot be applied \cite{he2016isoparametric}.

Efforts to extend isoparametric hypersurfaces through the introduction of anisotropic surface energy began somewhat oblivious to the Finslerian formalism, mainly pushed by the interest of other fields, such as analysis and material sciences; e.g. \cite{ge2012anisotropic} and its references. The explicit definition of an isoparametric function in a Finsler manifold first appeared in 2016 \cite{he2016isoparametric}. The authors of the article, He et al., considered similar equations to \ref{transnRiem} and \ref{isoRiem} for the nonlinear gradient and Laplacian operators but had to relax a bit the smoothness conditions of $f$, $a$, and $b$ to fit the smoothness of the Finsler metric.

A Finsler transnormal function is isoparametric precisely when each regular level set has constant nonlinear mean curvature \cite[Theorem~1.1]{he2016isoparametric}, which generalizes the result by Cartan for Riemannian manifolds \cite{cartan1938familles}. Besides, on a connected and forward complete Finsler manifold, the regular level sets of a transnormal function are parallel along the direction of the gradient \cite[Proposition~4.1]{he2016isoparametric}. However, the level sets of a Finsler transnormal function need not be equidistant \cite{alexandrino2019finsler, chen2022transnormal}, as opposed to the Riemannian case \cite{wang1987isoparametric}. In particular, an extension of Huygens' principle for wavefronts in Finsler spaces is a rather generic version \cite{dehkordi2019huygens}.

In Minkowski spaces (the Finslerian analogue to Euclidean spaces) with either the Busemann-Hausdorff or the Holmes-Thompson volume form, all the hyperplanes, the Minkowski hyperspheres, and the dual Minkowski cylinders are isoparametric hypersurfaces \cite[Theorem~1.2]{he2016isoparametric}. Actually, they are the only ones \cite[Theorems~5.1--5.5]{he2016isoparametric}, hence classifying them. In a different way, the preceding work on anisotropic isoparametric hypersurfaces in Euclidean spaces also achieved this classification \cite[Theorem~1.1]{ge2012anisotropic}. Despite the similarity to the Euclidean case, the isoparametric hypersurfaces of a Minkowski space are much more abundant, because even two dual Minkowski cylinders with the same radius may not be isometric to each other \cite[Example~2]{he2016isoparametric}.

For the past few years, the classification of isoparametric hypersurfaces (with respect to the Busemann-Hausdorff volume form) has progressed, from Randers-Minkowski spaces \cite{he2016isoparametric}, Funk spaces \cite{he2017isoparametric}, and Randers spheres \cite{xu2018isoparametric} to Randers space forms \cite{he2020classifications}. The key to such development is the equivalence between the isoparametric hypersurfaces of a Randers space and those of a Riemannian manifold. The correspondence is derived from the fact that the geodesics of strongly convex Randers metrics are exactly the solutions to Zermelo's navigation problem in Riemannian spaces subjected to mild wind \cite{bao2004zermelo}.

Generally, if $(M, F)$ is any Finsler space and $W$ is a vector field satisfying $F(-W) < 1$, then another Finsler metric emerges from navigation by
\begin{equation}
\label{Zermelodef}
F\left( \frac{v}{Z(v)} - W \right) = 1, \forall v \in TM \setminus 0.
\end{equation}
We call $Z$ a \emph{Zermelo metric} on $M$, and the pair $(F, W)$ its \emph{navigation data}. These metrics were introduced in \cite{shen2003finsler}, and their geodesics studied in \cite{huang2011geodesics}.

Under strong wind condition (i.e. $F(-W) > 1$), Equation~\ref{Zermelodef} originates two Finsler pseudo-metrics defined in the same conic subset of $TM$ over $M$, which are hence called \emph{conic} Finsler pseudo-metrics; one is positive-definite, and the other has Lorentz signature \cite{javaloyes2018some}. In critical wind (i.e. $F(-W) = 1$), only a positive-definite Finsler metric in a conic subset of $TM$ is obtained and simply referred to as a \emph{conic} Finsler metric; particularly, when $F$ is Riemannian, $Z$ is a Kropina metric \cite{yoshikawa2014kropina}. For an arbitrary wind condition, we say $Z$ is a \emph{conic Zermelo pseudo-metric}.

The extension of standard Finsler metrics to conic ones arises somewhat seamlessly in many interesting problems since the early stages of Finsler geometry. Even so, to the best of our knowledge, the first systematic account on the topic has been published for less than a decade \cite{javaloyes2014definition}. In \emph{conic} Minkowski spaces, the hyperplanes, the Minkowski hyperspheres, and the dual Minkowski cylinders are also isoparametric hypersurfaces \cite[Theorem 1.1]{he2022isoparametric}, but they are not all. For instance, a helicoid is an isoparametric hypersurface of a $3$-dimensional conic Minkowski $(\alpha,\beta)$-space \cite[Theorem 1.2]{he2022isoparametric}.

Roughly speaking, \emph{conic} Zermelo pseudo-metrics model the anisotropic mechanics of a system submitted to the influence of some vector field, as is obviously the case for Zermelo's navigation problem. Other applications include the spread of forest fires \cite{markvorsen2016finsler, javaloyes2021general}, the propagation of sound in a wind \cite{gibbons2011geometry}, the path of a seismic ray \cite{yajima2009finsler}, a model for a graphene sheet when an electric field is externally applied \cite{cvetivc2012graphene}, the deviation of a charged particle through an electromagnetic field \cite{cricsan2020finsler}, analogue models \cite{gibbons2009stationary, dehkordi2019huygens, javaloyes2021applications}, and general Finslerian spacetimes \cite{caponio2014wind, javaloyes2020definition, herrera2021stationary}.

The classification of isoparametric hypersurfaces extends to a Zermelo space $(M, Z)$ with the Busemann-Hausdorff volume, where the navigation data $(F, W)$ is composed of a Minkowski norm $F$ and a homothetic vector field $W$ \cite{dong2020isoparametric}. Similarly to Randers spaces, the crucial step is to show that $(M, Z)$ and $(M, F)$ have the same isoparametric hypersurfaces. Really, the equivalence of isoparametric hypersurfaces with respect to the Busemann-Hausdorff volume holds for a larger class of manifolds with navigation under the homothety condition of the wind vector field, as demonstrated in \cite{xu2020some} for any Zermelo metric, in \cite{xu2021isoparametric} for any conic Zermelo metric with Lorentz signature, and in \cite{he2022isoparametric} for a Kropina metric of constant flag curvature.

We prove that the isoparametric hypersurfaces of $(M, F)$ and $(M, Z)$ coincide for any conic Zermelo pseudo-metric -- with respect to some fixed volume form on $M$. Hence, our first theorem encompasses many previous results; e.g. \cite[Theorem 1.1]{he2020classifications}, \cite[Theorem 1.1]{dong2020isoparametric}, \cite[Theorem 1.5]{xu2020some}, \cite[Theorem 1.1]{xu2021isoparametric}, \cite[Theorem 1.3]{he2022isoparametric}. Nevertheless, the classification of isoparametric hypersurfaces is \emph{not} complete in many cases.

\begin{theorem}
\label{Zermeloisop}
Let $(M, Z, \nu)$ be a Finsler m space composed of a conic Zermelo pseudo-metric $Z$ with navigation data $(F, W)$ and a volume form $\nu$. Let $\varphi^W$ be the flow of the vector field $W$. Suppose $W$ is $F$-homothetic (i.e. $(\varphi_t^W)^{\ast} F = \e^{-\sigma t} F$, for some $\sigma \in \mathbb{R}$), and $\dv W$ is constant (i.e. $(\varphi_t^W)^{\ast} \nu = \e^{\tau t} \nu$, for some $\tau \in \mathbb{R}$). If $\rho$ and $\tilde\rho$ are $Z, F$-distance functions, respectively, and $L$ is a regular fiber for both with $\nabla^Z  \rho \vert_{L} = (\nabla^F \tilde\rho + W) \vert_L$, then there exist open sets $\mathcal{U} \subset (M, Z)$ and $\mathcal{\tilde U} \subset (M, F)$ with $ L\subset \mathcal{U} \cap \mathcal{\tilde U}$ such that $\rho\vert_{\mathcal{U}}$ is $(Z,\nu)$-isoparametric if and only if $\tilde\rho\vert_{\mathcal{\tilde U}}$ is $(F,\nu)$-isoparametric. In particular, $L\subset M$ is a $(Z,\nu)$-isoparametric hypersurface with respect to an admissible unit normal vector field $\xi$ if and only if it is $(F,\nu)$-isoparametric with respect to the unit normal vector field $\xi-W$.
\end{theorem}

The assumption over the divergence of $W$ is quite natural as it is satisfied for the Busemann-Hausdorff volume form associated with $Z$ when the vector field is $F$-homothetic; see Lemma~\ref{divBH}. In truth, the usual definition of the Busemann-Hausdorff volume form for a conic Zermelo metric deserves discussion because, under strong or critical wind, the domain of $Z$ at a point is a cone. For the present, we settle it by declaring $\nu^Z_{BH}$ to be the Busemann-Hausdorff volume form associated with $F$, once such equivalence holds for mild wind; e.g. \cite[Proposition~5.3]{shen2016introduction}.

The hypothesis over the gradients of the distance functions is also reasonable, although technical.  It concerns the idea that $\nabla^Z\rho$ and $\nabla^F\tilde\rho$ lie on the ``same side'' of the hypersurface $L$. This is particularly relevant for us because the two unitary normal vectors at a generic point of $L$ with respect to $Z$ are not multiples of each other -- due to the non-reversibility of the Zermelo metric, regardless if $F$ is reversible; see Remark~\ref{gradhypothesis} for details.

It now becomes clear that Theorem~\ref{Zermeloisop} generalizes those mentioned earlier, but the demonstration here is considerably different. While past cases compare the non-linear Laplacian for $(Z,\nu)$ to the Laplacian for the Riemannian metric determined by $\nabla^Z \rho$ with both evaluated at $\rho$, we will relate the non-linear Laplacians for $(Z, \nu)$ and $(F, \nu)$ of $\rho$ and $\tilde\rho$, respectively.

In the former technique, the constancy of the S-curvature is necessary, which is usually proved in local coordinates for the Busemann-Hausdorff volume. More specifically, the S-curvature turns up in the equation
\begin{equation}
\label{LaplaceScurv}
\triangle^Z\rho = \hat{\triangle}\rho - \mathrm{S} (\nabla^Z \rho) ,
\end{equation}
where $\hat{\triangle}$ denotes the Laplacian of $g^Z_{\nabla^Z \rho}$; see \cite[Lemma~14.1.2]{shen2001lectures} or  \cite[Proposition~3.1]{shen1997curvature} with opposite choice of sign for S. Our approach uses basic properties of the divergence -- once related the flows of the gradients of $\rho$ and $\tilde\rho$ via the flow of $W$; consequently, the proof becomes coordinate-free and valid for different volume forms, e.g. the Holmes-Thompson volume form associated with $F$ vide Lemma~\ref{divHT}. By thorough research, we find this path mostly unexplored, only lightly tread for Randers spheres \cite{chen2023isoparametric}. To pave way, we turn to positive-definite conic Finsler metrics and find the equations associating the non-linear Laplacians for $(Z, \nu)$ and $(F, \nu)$ of the same smooth function.

\begin{proposition}
\label{Laplacediv}
Let $(M, Z, \nu)$ be a conic Finsler m space, where $Z$ is a conic Zermelo metric with navigation data $(F, W)$ and $\nu$ some volume form. If $f $ is a smooth real function on some open $\mathcal{U} \subset (M, Z)$, then at regular admissible points
\begin{equation}
\label{Laplacefunc}
\begin{split}
&\frac{1}{Z(\nabla^Z f)} \left[ \triangle^Z f - \hs^Z f \left( \frac{\nabla^Z f}{Z(\nabla^Z f)}, \frac{\nabla^Z f}{Z(\nabla^Z f)} \right) \right] = \\
= &\frac{1}{F(\nabla^F f)} \left[ \triangle^F f - \hs^F f \left( \frac{\nabla^F f}{F(\nabla^F f)}, \frac{\nabla^F f}{F(\nabla^F f)} \right) \right] + \dv W .
\end{split}
\end{equation}
In particular, if $f$ is a $Z$-distance function, then
\begin{equation}
\label{Laplacedistfunc}
\triangle^Z f = \frac{1}{F(\nabla^F f)} \left( \triangle^F f - \hs^F f \left( \frac{\nabla^F f}{F(\nabla^F f)}, \frac{\nabla^F f}{F(\nabla^F f)} \right) \right) + \dv W ;
\end{equation}
if, additionally, $W$ is tangent to the fibers of $f$, then
\begin{equation}
\label{Laplacedist}
\triangle^Z f = \triangle^F f + \dv W.
\end{equation}
\end{proposition}

From Equation~\ref{Laplacefunc}, the function $f$ is $(Z, \nu)$-isoparametric if and only if it is $(F, \nu)$-isoparametric when $W$ and $\dv W$ are $f$-projectable, which adapts Theorem~\ref{Zermeloisop} to any smooth function under more geometric hypotheses on $W$ and $f$; see Corollary~\ref{Zermeloisopproj}. Next, Equation~\ref{Laplacedistfunc} offers an expression for the nonlinear mean curvature on hypersurfaces of $(M, Z, \nu)$ in terms of the navigation data because this curvature is, simply put, the nonlinear Laplacian of a distance function; vide \cite[Proposition~14.2]{shen1997curvature}. By further simplification, Equation~\ref{Laplacedist} resembles Equation~\ref{LaplaceScurv}, but the essential difference is that $\triangle^F$ is still a nonlinear operator without an explicit expression in general. Nonetheless, we investigate mean curvatures.

\begin{corollary}
\label{Zermelomean}
Let $(M, Z)$ be a conic Zermelo manifold with navigation data $(F, W)$ and let $\nu_{BH}$ be the Busemann-Hausdorff volume form associated with $Z$. If $L \subset M$ is a hypersurface and $\xi$ is an admissible unit normal vector field along $L$ with respect to $Z$, then $\xi-W$ is a unit normal vector field along $L$ with respect to $F$, and the mean curvature for $(M, Z, \nu_{BH})$ of $L$ in the direction of $\xi$ is
\begin{equation}
\label{Zermelomeaneq}
\Pi^Z_{\xi} = \Pi^F_{\xi - W} + \dv W,
\end{equation}
where $\Pi^F_{\xi - W}$ is the mean curvature for $(M, F, \nu_{BH})$ of $L$ in the direction of $\xi-W$.
\end{corollary}

By use of Equation~\ref{LaplaceScurv}, if $\rho$ is a distance function on an open subset $\mathcal{U} \subset (M, Z)$ such that $L = \rho^{-1}(c)$ and $\xi = \nabla^R \rho \vert_L$, then the Riemannian mean curvature $\hat{\Pi}$ with respect to $g^Z_{\nabla^Z \rho}$ satisfies
\begin{equation}
\label{ZermelomeanS}
\Pi^Z_{\xi} = \hat{\Pi}_{\xi} - \mathrm{S} (\xi) ;
\end{equation}
refer to \cite[Corollary~14.3.2]{shen2001lectures} or \cite[Proposition~14.1]{shen1997curvature} for the opposite sign of S. Comparably to Equation~\ref{ZermelomeanS}, Equation~\ref{Zermelomeaneq} has an interesting computative edge -- the divergence depends only on the volume form; see Example~\ref{nonlinearmeanex} for illustration. Our cost is the restriction to the Busemann-Hausdorff volume, while Equations~\ref{LaplaceScurv} and ~\ref{ZermelomeanS} hold for any regular Finsler m space. 

Whereas the nonlinear mean curvature arises from the variational problem for the volume induced to hypersurfaces by the ambient volume form and a unit normal vector field, the variation of the volume associated with the induced metric is distinct, as the volumes themselves need not coincide. The latter defines, for each volume notion, a linear mean curvature on any immersed submanifold $L^m \subset (M^n, F)$; e.g. \cite[Section~9.1]{shen2016introduction}. The linear mean curvature with respect to the Busemann-Hausdorff volume has been the first introduced \cite{shen1998finsler}; it is the only one we examine here, and what we mean whenever we say \emph{linear mean curvature}. The immersed submanifolds for which this curvature vanishes are still called \emph{BH-minimal} for clarity. Finally, we concentrate on expressions for the linear mean curvature when $F$ is Riemannian and $F(-W) < 1$, hence $(M, Z)$ is a Randers space. In this case, we replace $F$ with the corresponding Riemannian metric $h$ and the notation for the Randers metric from $Z$ to $R$.

\begin{theorem}
\label{Randerslinearmean}
Let $(M, R)$ be a Randers manifold with navigation data $(h, W)$ and let $\nu_{BH}$ be the Bussmann-Hausdorff volume form for $R$. Suppose $L^m \subset M^n$ is an immersed submanifold. If $\mathbf{n}$ is a unit normal vector field along $L$ with respect to $h$, then the linear mean curvature of $L$ in the direction of $\mathbf{n}$ is
\begin{equation}
\label{Randerslinearmeaneq}
\mathcal{H}(\mathbf{n}) = \Pi^{h}_{\mathbf{n}} + \mathcal{B}(\mathbf{n}) ,
\end{equation}
where $\Pi^h_{\mathbf{n}}$ is the Riemannian mean curvature of $L$ in the direction of $\mathbf{n}$, 
\begin{equation*}
\mathcal{B}(\mathbf{n}) := \frac{m h(W, \mathbf{n}) h(D_{\mathbf{n}} W, \mathbf{n})}{1 - h(W, \mathbf{n})^2} ,
\end{equation*}
and $D$ denotes the Levi-Civita connection for $h$. In particular, $\mathcal H(\mathbf{n})= \Pi^h_{\mathbf{n}} $ when the vector field $W$ is tangent to $L$ or $h$-Killing.
\end{theorem}

Once more, Equation~\ref{Randerslinearmeaneq} has meaningful computational value; cf. \cite[Equation~57]{shen1998finsler}, \cite[Section~3]{wu2007local}, \cite[Sections~2 and~3]{balan2007bh}, \cite[Equation~3.36]{dasilva2011minimal}, \cite[Section~3]{cui2014minimal}, \cite[Section~2]{cui2017nontrivial}. Its geometric nature made it easier to construct nontrivial BH-minimal surfaces in Example~\ref{linearmeanex}; cf. \cite{souza2003minimal, souza2004bernstein, wu2007local, dasilva2011minimal, dasilva2014helicoidal, cui2014minimal, cui2017nontrivial}. Furthermore, putting together Equations~\ref{Zermelomeaneq} and~\ref{Randerslinearmeaneq}, we find a necessary and sufficient condition for the constancy of the linear mean curvature on the fibers of a transnormal function to characterize an $(R, \nu_{BH})$-isoparametric function; see Corollary~\ref{Randersisopproj}. It is also simple to create a Randers manifold for which such characterization fails, as per Example~\ref{linearmeancounterex}. So the linear mean curvature does not seem like a good candidate to extend the definition of isoparametric hypersurfaces to \emph{isoparametric submanifolds}. We wonder which notion could be.
\section{Preliminaries}

Let $M$ be a smooth $n$-dimensional connected manifold, and $\pi: TM \to M$ its tangent bundle. Suppose an open subset $A \subset TM$ is \emph{conic} (i.e. $rv \in A$, $\forall v \in A, r>0$) and $A_p := A \cap T_p M$ is connected for each $p \in M$. A nonnegative continuous function $F: A \to [0, \infty)$ is a \emph{conic Finsler metric} if $\pi(A) = M$, $F$ is smooth on $A \setminus 0$, it is positive homogeneous of degree one (i.e. $ F(rv) = rF(v)$, $\forall v \in A, r>0$), and for each $v \in A \setminus 0$ the symmetric bilinear form
\begin{equation*}
g_{v}(u_1, u_2) :=  \left.\frac{1}{2}\frac{\partial^{2}}{\partial t_2 \partial t_1} F^{2}(v + t_1 u_1 + t_2 u_2)\right\vert_{t_1=t_2=0}, \; \forall u_1, u_2\in T_{\pi(v)}M ,
\end{equation*}
is positive-definite. Clearly, $F$ is a Finsler metric in the usual sense when $A=TM$; in this case, we may refer to $F$ as convex, regular, or just a Finsler metric.

From the definition, $g$ is a Riemannian metric on the pullback bundle $\pi^{\ast} TM$ over $A \setminus 0$, called the \emph{fundamental tensor}. By homogeneity, $F^2(v) = g_v(v, v)$, $\forall v \in A \setminus 0$; or more generally,
\begin{equation*}
\frac{1}{2} \mathrm{d} \left( F^2 (v) \right)(u) = g_v(v, u), \, \forall v \in A \setminus 0 , u \in T_{\pi(v)}M .
\end{equation*}

The \emph{Cartan tensor} is another symmetric section on $\pi^{\ast}TM$, given by
\begin{equation*}
C_v(u_1, u_2, u_3) :=  \left.\frac{1}{4}\frac{\partial^{3}}{\partial t_3 \partial t_2 \partial t_1} F^{2}\left(v + t_1 u_1 + t_2 u_2 + t_3 u_3 \right)\right\vert_{t_1=t_2=t_3=0} ,
\end{equation*}
$\forall v \in A \setminus 0$, and $\forall u_1, u_2, u_3 \in T_{\pi(v)}M$. Homogeneity now implies $C_v(v, \cdot, \cdot) = 0$.

A nonzero vector $v \in A_p$ is \emph{orthogonal} to $u \in T_p M$ when $g_v(v, u) = 0$. Further, $v \in A_p \setminus 0$ is orthogonal to a submanifold $L \subset M$ if $p \in L$ and $v$ is orthogonal to all vectors in $T_p L$. Similarly, we may say $v$ is orthogonal to a submanifold $L\subset T_p M$ through the identification of $T_v L$ as a subspace of $T_p M$. In case $L \subset M$ is a hypersurface, we call $L$ \emph{admissible} when there exists an orthonormal smooth vector field $\xi$ along $L$ (particularly, $\xi(p) \in A _p \setminus 0$, $\forall p\in L$).

The \emph{indicatrix} of $F$ is the smooth connected hypersurface $\Sigma := F^{-1}(1)$ embedded in $A \subset TM$. Fixed a point $p \in M$, the indicatrix of $F$ at $p$ is $\Sigma_p := \Sigma \cap T_pM \subset A_p$. Each unit vector $v \in A_p$ is orthogonal to $\Sigma_p$, so $T_v \Sigma_p = \{ u \in T_pM : g_v (v, u) = 0 \}$.

A vector field $V$ over an open subset $\mathcal{U} \subset (M, F)$ is \emph{admissible} when $V(p) \in A_p \setminus 0$ for all $p \in \mathcal{U}$. For example, the orthonormal vector field $\xi$ along an admissible hypersurface $L \subset M$ is admissible. For regular Finsler metrics, the admissible vector fields are simply those nonvanishing.

The \emph{Chern connection} for $(M,F)$ can be viewed as a family of affine connections $\{ \nabla^V : V \text{ admissible vector field}\}$ that is \emph{torsion-free}:
\begin{equation*}
\nabla^V_X Y - \nabla^V_Y X = [X,Y] ;
\end{equation*}
and \emph{almost metric compatible}:
\begin{equation*}
X g_V(Y,  Z) = g_V(\nabla^V_X Y, Z) + g_V(Y, \nabla^V_X Z) + 2C_V(\nabla^V_X V, Y, Z) ;
\end{equation*}
here $X,Y,Z$ are smooth vector fields. By \cite[Proposition~2.6]{javaloyes2014chern}, $\nabla^V$ depends pointwise on $V$; i.e. if $v \in A_p \setminus 0$ and $V_1, V_2$ are admissible vector fields with $V_1(p) = V_2(p) = v$, then $(\nabla^{V_1}_{X} Y)(p) = (\nabla^{V_2}_{X} Y)(p)$ for all smooth vector fields $X,Y$. Of course, $\nabla^V_X Y$ also depends pointwise on $X$ for each $V$ fixed.

If $\gamma: I \to M$ is a smooth curve, and $V$ is an admissible vector field along $\gamma$, then $\nabla^V_{\gamma^{\prime}}$ provides a covariant derivative along $\gamma$ with reference $V$; e.g. \cite[Remark~2.7]{javaloyes2014chern}. A smooth curve $\gamma: I \to M$ is \emph{admissible} when $\gamma^{\prime}$ is an admissible vector field along $\gamma$, and an admissible curve $\gamma$ is a \emph{(nonconstant) geodesic} if $\nabla_{\gamma^{\prime}}^{\gamma^{\prime}} \gamma^{\prime} = 0$; by almost metric compatibility, $\gamma$ has constant speed. For each $v \in A_p \setminus 0$, there is a unique (nonconstant) geodesic $\gamma_v: I_v \to M$ such that $\gamma_v(0) = p$ and $\gamma_v^{\prime}(0) = v$.

A smooth function $f: \mathcal{U} \subset M \to \mathbb{R}$ is \emph{admissible} if there is an admissible vector field $\nabla f$ defined implicitly on $\{ p \in \mathcal{U} : \mathrm{d} f_p \neq 0 \}$ by
\begin{equation*}
g_{\nabla f_p} (\nabla f_p, u) = \mathrm{d} f_p (u) , \forall u \in T_pM ;
\end{equation*}
for points where $\mathrm{d} f_p = 0$, set $\nabla f(p) = 0$. The resulting vector field $\nabla f$ is the \emph{gradient} of $f$. Although the gradient is not linear, it is still orthogonal to the regular level sets of $f$. When $A = TM$, all smooth functions are admissible, and the gradient is well-defined for any smooth function.

A nonconstant admissible function $f: \mathcal{U} \subset M \to \mathbb{R}$ is \emph{transnormal} when $F^2 (\nabla f) = b \circ f$ for some real continuous function $b$ on $f(\mathcal{U})$ which is smooth on the interior of $f (\{ p \in \mathcal{U} : \mathrm{d} f_p \neq 0 \})$. A locally Lipschitz  function $\rho: \mathcal{U} \subset M \to [0, \infty)$ is a \emph{distance function} if it is admissible almost everywhere and $F (\nabla \rho) = 1$ at all regular admissible points; in particular, the restriction of $\rho$ to the set of regular admissible points is a transnormal function.

\begin{lemma}
\label{transnormalintegralcurve}
Let $(M,F)$ be a conic Finsler space. If $f$ is a transnormal function and $\gamma$ is an integral curve of $\nabla f$, then $\gamma$ is a reparametrization of a geodesic. For a distance function, $\gamma$ is a piecewise geodesic of unit speed.
\end{lemma}

The demonstration for this result is easily adapted from \cite[Proposition~4.4(a)]{alexandrino2019finsler}, and hence will be omitted.

\begin{lemma}
\label{existofdistfunc}
Let $(M,F)$ be a conic Finsler space, and let $L \subset M$ be an admissible hypersurface. If $\xi$ is an orthonormal smooth vector field along $L$, then for each $p \in L$ there exists an open neighborhood $\mathcal{U} \subset (M, F)$ around $p$, and a distance function $\rho: \mathcal{U} \to \mathbb{R}$ such that $\nabla \rho\vert_{\mathcal{U} \cap L} = \xi\vert_{\mathcal{U} \cap L}$. Furthermore, the function $\rho$ is unique up to a constant.
\end{lemma}
\begin{proof}
Fixed $p \in L$, there is an open neighborhood $\mathcal{\tilde U} \subset (M, F)$ around $p$ and $\varepsilon > 0$ such that
\begin{equation*}
H(q,t):= \gamma_{\xi(q)}(t)\, , \; \forall (q, t) \in (\mathcal{\tilde U} \cap L) \times (-\varepsilon, \varepsilon),
\end{equation*}
is a well-defined smooth map. Moreover, $\mathrm{d} H_{(p,0)}$ is an linear isomorphism. Then, by the Inverse Function Theorem and reducing $\mathcal{\tilde U}$ if necessary, we can assume $H$ is a diffeomorphism onto $\mathcal{U} := \im H$.

Consider the function $\rho: \mathcal{U} \to \mathbb{R}$ defined by $\rho(\gamma_{\xi(q)} (t)) = t$. Using almost metric compatibility of the Chern connection and proceeding analogously to the Riemannian case, we conclude that  $\nabla \rho_{\gamma_{\xi(q)} (t)} = \gamma^{\prime}_{\xi(q)} (t)$. Therefore, $\rho$ is a distance function with $\nabla \rho\vert_{\mathcal{U} \cap L} = \xi\vert_{\mathcal{U} \cap L}$. 
  
Lastly, uniqueness follows from Lemma~\ref{transnormalintegralcurve}.
\end{proof}

The proof of Lemma~\ref{existofdistfunc} is a standard one but it needs to be carefully formulated because there are subtleties to the general properties of conic Finsler spaces; find detailed discussion in \cite[Section~3]{javaloyes2014definition}.

Let $f: \mathcal{U} \subset M \to \mathbb{R}$ be an admissible function. The \emph{Hessian} of $f$ is a symmetric bilinear map defined on $\{ p \in \mathcal{U} : \mathrm{d} f_p \neq 0 \}$ by
\begin{equation*}
\hs f (u,v) := g_{\nabla f_p} (\nabla^{\nabla f_p}_u \nabla f_p, v) , \, \forall u,v \in T_p M.
\end{equation*}
For smooth vector fields $X, Y$ on $\{ p \in \mathcal{U} : \mathrm{d} f_p \neq 0 \}$, we have
\begin{equation*}
\hs f (X,Y) = X Y (f) - \mathrm{d} f (\nabla^{\nabla f}_X Y) .
\end{equation*}

Now let $f: \mathcal{U} \subset M \to \mathbb{R}$ be any continuous function. Another real function $\phi$ defined over the same open $\mathcal{U}$ is \emph{projectable} with respect to $f$ -- or \emph{$f$-projectable} for short -- if $\phi$ is constant on the fibers of $f$. A smooth vector field $X$ on $\mathcal{U}$ is \emph{$f$-projectable} if its flow takes any fiber of $f$ into a fiber of $f$. When $f$ is differentiable, $X$ is $f$-projetable if only if the function $\mathrm{d} f (X)$ is $f$-projetable. In this language, a nonconstant admissible function $f$ is transnormal if only if $F(\nabla f)$ is $f$-projectable.

\begin{lemma}
\label{Hessproj}
Let $(M, F)$ be a conic Finsler space. If $f: \mathcal{U} \subset M \to \mathbb{R}$ is a transnormal function with $F^2(\nabla f) = b \circ f$, then
\begin{equation*}
\hs f(\nabla f, \nabla f) = \frac{1}{2} b^{\prime}(f) b(f) .
\end{equation*}
In particular, $\hs f (\nabla f, \nabla f)$ is $f$-projectable.
\end{lemma}

The argument in \cite[Lemma~6.6]{alexandrino2019finsler} holds with minor adjustments, so we omit the demonstration again.

Suppose $M$ is endowed with a volume form $\nu$. Accurately, $\nu$ is a family of non-degenerate top-dimensional forms $\nu_i$ on coordinate charts of an atlas $\{(x_i, \mathcal{U}_i)\}$ for $M$ such that $\nu_i = \nu_j$ on $\mathcal{U}_i \cap \mathcal{U}_j$ whenever the intersection is non-empty; that is, if $\nu_i = \mu_i (x) dx^1 \cdots dx^n$ and $\nu_j = \mu_j (\tilde x) d\tilde x^1 \cdots d\tilde x^n$, then $\mu_i (x) = \left\vert \det \left( \frac{\partial \tilde x^k}{\partial x^l} \right) \right\vert \mu_j (\tilde x)$.

The divergence of a smooth vector field $X$ in $(M, \nu)$ is the function defined by
\begin{equation*}
(\dv X) \nu := \mathrm{d} (i_X \nu) = \left. \frac{d}{d t} \varphi^{\ast}_t \nu \right\vert_{t = 0} = \mathcal{L}_X \nu,
\end{equation*}
where $i_X$ denotes the inclusion of $X$ at the first entry, $\varphi$ is the flow of $X$, and $\mathcal{L}$ denotes the Lie derivative on $M$. If $X, Y$ are smooth vector fields and $f$ is a smooth function, then the following properties hold:
\begin{subequations}
\begin{equation*}
\dv (X + Y) = \dv X + \dv Y ;
\end{equation*}
\begin{equation*}
\dv (f X) = f \dv X + \mathrm{d} f (X) .
\end{equation*}
\end{subequations}

Given a conic Finsler metric $F$ and a volume form $\nu$ on $M$, the triple $(M, F, \nu)$ is called a \emph{conic Finsler m space}.

\begin{lemma}
\label{divhypersurface}
Let $(M, F, \nu)$ be a conic Finsler m space, and let $L \subset M$ be an admissible hypersurface. If $X_1, X_2$ are admissible unit vector fields on some neighborhood of $L$ such that $X_1\vert_L = X_2\vert_L$ is orthogonal to $L$, then $(\dv X_1)|_L = (\dv X_2)|_L$.
\end{lemma}
\begin{proof}
For each $p \in L$, take a local frame $\{E_1, \ldots, E_n\}$ for $TM$ around $p$ such that $\{E_1\vert_L, \ldots, E_{n-1}\vert_L\}$ is a local frame for $TL$ and $E_n := X_2$.

On one hand, since $X_1\vert_L = X_2\vert_L$, we have that $(\nabla^{X_1}_{X_1} E_i)(p) = (\nabla^{X_2}_{X_2} E_i)(p)$ and $(\nabla^{X_1}_{E_i} X_1)(p) = (\nabla^{X_2}_{E_i} X_2)(p)$ for $i = 1, \ldots, n-1$. Hence, $[X_1, E_i](p) = [X_2, E_i](p)$ for $i = 1, \ldots, n-1$. On the other hand, $g_{X_j}(X_j, X_j) = 1$ for $j = 1, 2$ implies that $g_{X_j}(X_j, \nabla^{X_j}_{X_k} X_j) = 0$ for $j, k = 1, 2$. So $(\nabla^{X_j}_{X_k} X_j)(p)$ is tangent to $L$. Using $X_1\vert_L = X_2\vert_L$ again, we conclude $[X_1, X_2](p)$ is tangent to $L$. Consequently,
\begin{equation*}
\nu_p (E_1, \ldots, E_{n-1}, [X_1, E_n]) = \nu_p (E_1, \ldots, E_{n-1}, [X_2, E_n]) = 0 .
\end{equation*}
 
Now
\begin{align*}
(\mathcal{L}_{X_1} \nu)_p (E_1, \ldots, E_n) &= X_1 (\nu (E_1, \ldots, E_n))(p) - \sum_{i=1}^n \nu_p (E_1, \ldots, [X_1, E_i], \ldots, E_n) \\
&= X_2 (\nu (E_1, \ldots, E_n))(p) - \sum_{i=1}^n \nu_p (E_1, \ldots, [X_2, E_i], \ldots, E_n)\\
&= (\mathcal{L}_{X_2} \nu)_p (E_1, \ldots, E_n) .
\end{align*}
Therefore, $(\mathcal{L}_{X_1} \nu)_p = (\mathcal{L}_{X_2} \nu)_p$, whence $(\dv X_1)(p) = (\dv X_2)(p)$.
\end{proof}

For a Riemannian manifold $(M, h)$ there is only one natural choice for the volume form; namely,
\begin{equation*}
\nu^h := \sqrt{\det h_{ij} (x)} dx^1 \cdots dx^n ,
\end{equation*}
where $h_{ij} (x)$ are the components of the metric with respect to the local basis induced by coordinates. In general, there is no preferred choice; the most important volume forms associated with a Finsler metric are the Busemann-Hausdorff and the Holmes-Thompson, which we introduce next -- based on the direct approach of \cite[Section~2.2]{shen2001lectures}. When the Finsler space is conic, a volume form will be fixed subjectively as needed.

Assume $(M, F)$ is an oriented Finsler space. For each $p \in M$, take an arbitrary oriented basis $\{e_i\}_{i=1}^{n}$ of $T_pM$ and the dual basis $\{\theta^i\}_{i=1}^{n}$ for $T^{\ast}_pM$.

Let $B^n (p) := \{(y^i) \in \mathbb{R}^n : F(y^i e_i) < 1 \}$, $\mathbb{B}^n := \{(y^i) \in \mathbb{R}^n : \sum (y^i)^2 < 1 \}$, and let $\vl (\Omega)$ denote the Euclidean volume of a bounded open subset $\Omega \subset \mathbb{R}^n$. Consider the function
\begin{equation}
\label{BH}
\mu_{BH} (p) := \frac{\vl (\mathbb{B}^n)}{\vl (B^n (p))} .
\end{equation}
The \emph{Busemann-Hausdorff volume form} for $(M, F)$ is
\begin{equation*}
(\nu_{BH})_p := \mu_{BH} (p) \theta^1 \wedge \cdots \wedge \theta^n .
\end{equation*}
The measure it generates on $M$ is the Hausdorff measure of the distance function determined by $F$ when the metric is reversible (i.e. $F(-v) = F(v)$, $\forall v \in TM$).

For each $y = y^ie_i \in T_pM \setminus \{0\}$, the components of the fundamental tensor $g_{ij} (y) := g_{y} (e_i, e_j)$ are smooth functions on $\mathbb{R}^n \setminus \{0\}$. Consider now the function
\begin{equation}
\label{HT}
\mu_{HT} (p) := \frac{\int_{B^n(p)} \det( g_{ij} (y) ) dy^1 \cdots dy^n}{\vl (\mathbb{B}^n)} .
\end{equation}
The \emph{Holmes-Thompson volume form} associated with $(M, F)$ is
\begin{equation*}
(\nu_{HT})_p := \mu_{HT} (p) \theta^1 \wedge \cdots \wedge \theta^n .
\end{equation*}
It arises from the volume form $\hat{\nu} := (-1)^{\frac{n(n+1)}{2}} \frac{1}{n!} (\mathrm{d} \omega)^n$ on the unit ball bundle of $M$, where $\omega:= g_{ij}(y)y^j dx^i$ is the Hilbert form.

Allow $(M, F, \nu)$ to be any conic Finsler m space. If $f: \mathcal{U} \subset M \to \mathbb{R}$ is an admissible function, then the \emph{(nonlinear) Laplacian} of $f$ in $(M, F, \nu)$ -- or simply \emph{$(F, \nu)$-Laplacian} -- is defined on $\{p \in \mathcal{U} : \mathrm{d} f_p \neq 0 \}$ by
\begin{equation*}
\triangle f := \dv (\nabla f) .
\end{equation*}
Finally, the function $f$ is \emph{isoparametric} in $(M, F, \nu)$ -- or \emph{$(F, \nu)$-isoparametric} -- when there are real continuous functions $a, b$ on $f(\{p \in \mathcal{U} : \mathrm{d} f_p \neq 0 \})$ with $b$ smooth over the open domain such that
\begin{align*}
F^2 (\nabla f) &= b \circ f ; \\
\triangle f &= a \circ f .
\end{align*}
In this case, $f$ is transnormal as per the definition above by continuity.

Descriptively, $f$ is an $(F,\nu)$-isoparametric function if only if $F (\nabla f)$ and $\triangle f$ are $f$-projectable. In particular, a transnormal function $f$ is $(F, \nu)$-isoparametric if and only if $\triangle f$ is $f$-projectable.

A hypersurface $L\subset M$ is an \emph{$(F, \nu)$-isoparametric hypersurface} with respect to a unit normal vector field $\xi$ along $L$ if there is an $(F, \nu)$-isoparametric function $f$ for which $L$ is a regular level set and $\left. \frac{\nabla f}{F(\nabla f)} \right\vert_L = \xi$.
\section{Zermelo Metrics}

Suppose $M$ is a manifold as before, $F: TM \to [0, \infty)$ is a Finsler metric, and $W$ is a smooth vector field on $M$. Denote the indicatrix of $F$ by $\Sigma^F$ and its fundamental tensor by $g^F$. Let $\Sigma := \{u + W : u \in \Sigma^F , g^F_u(u, -W) < 1\}$ and $A \setminus 0 := \{rv : v \in \Sigma , r>0\}$. Define a function $Z: A \to [0, \infty)$ by
\begin{equation*}
Z(v) = r > 0 \text{ if and only if } r^{-1}v \in \Sigma, \text{ and } Z(v) = 0 \text{ if and only if } v = 0 .
\end{equation*}
The positive values of $Z$ are expressed implicitly by
\begin{equation*}
F\left( \frac{v}{Z(v)} - W \right) = 1 , \forall v \in A \setminus 0 .
\end{equation*}
Therefrom, we obtain a well-defined positive-definite conic Finsler metric $Z$, whose indicatrix is $\Sigma$. We call $Z$ a \emph{conic Zermelo metric} with \emph{navigation data} $(F, W)$.

The Zermelo metric is a regular Finsler metric if and only if $F(-W) < 1$. The subset $\tilde \Sigma := \{u + W : u \in \Sigma^F , g^F_u(u, -W) > 1\}$ is non-empty if and only if $F(-W) > 1$. In this case, the above construction for $\tilde \Sigma$ yields the same conic subset $A \subset TM$ and a conic Zermelo pseudo-metric $Z$ with Lorentz signature; see details in \cite{javaloyes2018some}. We assume the choice of $\Sigma$ (or positive-definite conic metrics) unless stated otherwise. For a Riemannian norm $F = \sqrt{h}$, the metric $Z$ is of Randers type, denoted distinctively by $R$, and the pair $(h, W)$ is the navigation data.

\begin{remark}
Once we work with two different metrics -- $F$ and $Z$ --, we use prefix and superscript as indicators of the metric to which a quantity refers. For instance, the fundamental tensor of $Z$ is $g^Z$, the gradient of an admissible function $f$ with respect to the Zermelo metric is the $Z$-gradient of $f$, denoted $\nabla^Z f$, and so on.
\end{remark}

\begin{lemma}
\label{Zermelofundtensor}
Let $Z$ be a conic Zermelo pseudo-metric on $M$ with navigation data $(F, W)$. For all $v \in A \setminus 0$ and $u \in T_{\pi(v)}M$,
\begin{equation*}
g^Z_v (v, u) = k(v) g^F_{v - Z(v)W} (v - Z(v)W, u),
\end{equation*}
where the function $k$ is nowhere vanishing on $A \setminus 0$ and positive homogeneous of degree zero. In particular, $v$ is $Z$-orthogonal to $u$ if and only if $v - Z(v)W$ is $F$-orthogonal to $u$.
\end{lemma}
\begin{proof}
For each $v \in A \setminus 0$, $Z(v) = r$ if and only if $F(v - Z(v)W) = r$, with $r>0$.

Let $p := \pi(v)$ and consider the following linear functionals on $T_p M$:
\begin{equation*}
\ell^Z_v := g^Z_v(v, \cdot) \, , \; \ell^F_{v - Z(v)W} := g^F_{v - Z(v)W}(v - Z(v)W, \cdot) \, .
\end{equation*}
Their kernels are $T_v(r\Sigma^Z_p)$ and $T_{v - Z(v)W}(r\Sigma^F_p)$, which coincide as vector subspaces of $T_pM$. So there exists a real number $k(v)$ such that $\ell^Z_v = k(v) \, \ell^F_{v - Z(v)W}$.

By definition,
\begin{equation*}
k(v) = \frac{g^Z_v(v, v - Z(v)W)}{g^F_{v - Z(v)W}(v - Z(v)W, v - Z(v)W)} \, .
\end{equation*}
Since $F(v - Z(v)W) = Z(v)$,
\begin{equation*}
k(v) = \frac{g^Z_v(v, v) + Z(v)g^Z_v(v, -W)}{Z^2(v)} = 1 + \frac{g^Z_v(v, -W)}{Z(v)} \, .
\end{equation*}
In particular, $k$ is positive homogeneous of degree zero.

Now,
\begin{equation*}
g^Z_v(v, -W) = \left[ 1 + \frac{g^Z_v(v, -W)}{Z(v)} \right] g^F_{v - Z(v)W}(v - Z(v)W, -W) \, ,
\end{equation*}
whence
\begin{equation*}
g^F_{v - Z(v)W}(v - Z(v)W, -W) = g^Z_v(v, -W) \left[ 1 - g^F_{\frac{v}{Z(v)} - W}\left(\frac{v}{Z(v)} - W, -W\right) \right] \, .
\end{equation*}
Since $g^F_{\frac{v}{Z(v)} - W}\left(\frac{v}{Z(v)} - W, -W\right) < 1$ or $> 1$,
\begin{equation*}
g^Z_v(v, -W) = \frac{Z(v) g^F_{\frac{v}{Z(v)} - W}\left(\frac{v}{Z(v)} - W, -W\right)}{1 - g^F_{\frac{v}{Z(v)} - W}\left(\frac{v}{Z(v)} - W, -W\right)} \, .
\end{equation*}

Therefore,
\begin{equation*}
k(v) = 1 + \frac{g^F_{\frac{v}{Z(v)} - W}\left(\frac{v}{Z(v)} - W, -W\right)}{1 - g^F_{\frac{v}{Z(v)} - W}\left(\frac{v}{Z(v)} - W, -W\right)} = \frac{1}{1 - g^F_{\frac{v}{Z(v)} - W}\left(\frac{v}{Z(v)} - W, -W\right)} \, .
\end{equation*}
\end{proof}

This result generalizes \cite[Lemma~3.1]{alexandrino2019finsler} to any conic Zermelo pseudo-metric. Notice that $k$ is a positive function on $A\setminus 0$ whenever the Zermelo metric is positive-definite. Otherwise, $Z$ is a pseudo-metric with Lorentz signature and the function $k$ is negative on $A \setminus 0$. Lastly, for any conic Finsler space $(M, F)$ and $p \in M$, the map $\ell: T_pM \cap A \setminus 0 \to T_pM^\ast$ given by the correspondence $v \mapsto \ell_v := g^F_v(v, \cdot)$ is called the \emph{Legendre transformation} at $p$. If the Finsler metric is regular, then $\ell$ is a bijection from $T_pM$ to $T_pM^\ast$. But in the conic case, $\ell$ is not surjective.

\begin{lemma}
\label{Zermelograd}
Let $(M, Z)$ be a conic Zermelo space with navigation data $(F, W)$. If $f: \mathcal{U} \subset M \to \mathbb{R}$ is a $Z$-admissible function, then
\begin{enumerate}
\item[(a)] $\frac{\nabla^Z f}{Z(\nabla^Z f)} = \frac{\nabla^F f}{F(\nabla^F f)} + W$ on $\{p \in \mathcal{U} : \mathrm{d} f_p \neq 0\}$; and
\item[(b)] $Z(\nabla^Z f) = F(\nabla^F f) + \mathrm{d} f(W)$.
\end{enumerate} 
\end{lemma}

The proof is analogous to the one of \cite[Lemma~3.2]{alexandrino2019finsler}. But it is essential that $g^Z$ is positive-definite.

\begin{corollary}
\label{ZFtransnormal}
Let $(M, Z)$ be a conic Zermelo space with navigation data $(F, W)$. For any $Z$-admissible function $f$, the following are equivalent:
\begin{enumerate}
\item $W$ is $f$-projetable;
\item $f$ is $Z$-transnormal if only if it is $F$-transnormal.
\end{enumerate}
\end{corollary}

When $Z$ is a regular Finsler metric on $M$ arising from $(F, W)$ via navigation, the Busemann-Hausdorff volume forms of $Z$ and $F$ are the same, i.e.

\begin{proposition}[Proposition~5.3,\cite{shen2016introduction}]
\label{ZermeloBH}
If $(M, Z)$ is a Zermelo space with navigation data $(F, W)$ and $F(-W) < 1$, then $\nu^Z_{BH} = \nu^F_{BH}$.
\end{proposition}

In the conic case, the definition of the Busemann-Hausdorff volume form must be revised, since not all vectors are admissible. But we avoid this discussion here. For our purposes, it suffices to declare that $\nu^Z_{BH}$ is $\nu^F_{BH}$, based on the above.

The last result of this section is an important step in the proof of Theorem~\ref{Randerslinearmean}.

\begin{lemma}
\label{Randersimmersion}
Let $(M, R)$ be a Randers manifold with navigation data $(h, W)$. If $L^m \subset M^n$ is an immersed submanifold, then the induced metric on $L$ is the Randers metric with navigation data $(\tilde h, W^{\top})$ given by
\begin{subequations}
\begin{equation*}
\tilde h := \lambda \, h\vert_{TL \times TL} \, ,
\end{equation*}
\begin{equation*}
W^{\top} := W - W^{\perp} \, ,
\end{equation*}
\end{subequations}
where $\lambda := \frac{1}{1 - h(W^{\perp}, W^{\perp})}$ and $W^{\perp}$ is the component of $W$ $h$-orthogonal to $L$.
\end{lemma}
\begin{proof}
If $v \in TL \setminus 0$, then
\begin{equation*}
1 = h \left( \frac{v}{R(v)} - W, \frac{v}{R(v)} - W \right) = h \left( \frac{v}{R(v)} - W^{\top}, \frac{v}{R(v)} - W^{\top} \right) + h \left( W^{\perp}, W^{\perp} \right) .
\end{equation*}
Since $h(W^{\perp}, W^{\perp}) \leq h(W,W) <1$,
\begin{equation*}
\frac{1}{1 - h(W^{\perp}, W^{\perp})} h \left( \frac{v}{R(v)} - W^{\top}, \frac{v}{R(v)} - W^{\top} \right) = 1 .
\end{equation*}
So the induced metric $R\vert_{TL}$ is the Randers metric with navigation data $(\tilde h, W^{\top})$.
\end{proof}

\begin{corollary}
\label{RandersimmersionBH}
Let $(M, R)$ be a Randers manifold with navigation data $(h, W)$ and let $\nu^h$ be the Riemannian volume form of $h$. If $L^m \subset M^n$ is an immersed submanifold, then the Busemann-Hausdorff volume form associated with the induced metric on $L$ is $\lambda^{\frac{m}{2}} (\nu^h)\vert_{TL}$, for $\lambda$ as in Lemma~\ref{Randersimmersion}.
\end{corollary}
\begin{proof}
By Proposition~\ref{ZermeloBH}, the wanted volume form is the Riemannian volume form associated with $\tilde h$. So Lemma~\ref{Randersimmersion} concludes the proof.
\end{proof}
\section{Theorem \ref{Zermeloisop}}

We first need a geometric description for the geodesics of a conic Zermelo pseudo-metric $Z$ in terms of its navigation data $(F,W)$; to be precise, when $W$ is an $F$-homothetic vector field we have the following theorem.

\begin{theorem}[Theorem~1.2, \cite{javaloyes2018some}]
\label{Zermelogeod}
Let $(M,Z)$ be a conic Zermelo pseudo-metric space with navigation data $(F,W)$. If $W$ is an $F$-homotetic vector field with coefficient $-\sigma$, then the unit speed geodesics of $Z$ are expressed on some neighborhood of $t=0$ by
\begin{equation*}
 \gamma(t) = \varphi^W_t( \tilde{\gamma}( s(t) ) ) \, ,
\end{equation*}
where $\varphi^W$ is the flow of $W$, $\tilde{\gamma}$ is a unit speed geodesic of $(M,F) $ and
\begin{equation*}
s(t) = \left\{ \begin{array}{cc} \frac{e^{\sigma t} - 1}{\sigma } , & \text{ if } \sigma \neq 0 ; \\ t , & \text{ if } \sigma = 0. \end{array} \right.
\end{equation*}
\end{theorem}

This result was initially proved for Randers metrics \cite{robles2007geodesics}. It was later extended to Zermelo metrics arising from a navigation problem under mild wind condition (i.e. $F(-W) < 1$) \cite{huang2011geodesics}, and generalized to the conic case for pseudo-metrics \cite{javaloyes2018some}.

\begin{proof}[Proof of Theorem~\ref{Zermeloisop}]
Given an $F$-distance function $\tilde \rho$ and a regular fiber $L$ of $\tilde \rho$, there is an open subset $\tilde{\mathcal{U}} \subset (M, F)$ with $L \subset \tilde{\mathcal{U}}$ such that $\tilde \rho$ is regular at all points of $\tilde{\mathcal{U}}$. Let $\tilde \varphi$ be the flow of $\nabla^F \tilde \rho$. For each $p \in L$, the curve $t \mapsto \tilde \varphi_t (p)$ in $\tilde{\mathcal{U}}$ is a unit speed geodesic of $(M, F)$ with initial velocity $\nabla^F \tilde \rho (p)$.

Similarly, there exists an open $\mathcal{U} \subset (M, Z)$ with $L \subset \mathcal{U}$ such that $\rho$ is regular on $\mathcal{U}$ and, if $\varphi$ is the flow of $\nabla^Z \rho$, then the curve $t \mapsto \varphi_t (p)$ in $\mathcal{U}$ is a unit speed geodesic of $Z$ with initial velocity $\nabla^Z \rho (p)$.

By Theorem~\ref{Zermelogeod}, the curve $t \mapsto \varphi^W_t \circ \tilde \varphi_{s(t)} (p)$ for $t$ small enough is another unit speed geodesic of $(M, Z)$ with initial velocity $\nabla^F \tilde \rho (p) + W(p) = \nabla^Z \rho (p)$. Since $\varphi^W$ is a local diffeomorphism, we may assume $\tilde \varphi_{s(t)} (p) \in \tilde{\mathcal{U}}$ if and only if $\varphi^W_t \circ \tilde \varphi_{s(t)} (p) \in \mathcal{U}$, taking smaller open subsets if necessary. So $\varphi_t (p) = \varphi^W_t \circ \tilde \varphi_{s(t)} (p)$ by uniqueness of geodesics.

For $\varphi_t(p) \in \mathcal{U}$, we have
\begin{equation*}
\nabla^Z \rho (\varphi_t (p)) = \frac{d}{dt} \varphi_t (p) = W(\tilde \varphi_{s(t)} (p)) + s^{\prime} (t) \mathrm{d} \varphi_t^W (\nabla^F \tilde \rho (\tilde \varphi_{s(t)} (p))) .
\end{equation*}
By divergence properties,
\begin{equation*}
\dv (\nabla^Z \rho) (\varphi_t (p)) = \dv W (\tilde \varphi_{s(t)} (p)) + s^{\prime\prime} (t) + s^{\prime}(t) \dv (\mathrm{d} \varphi_t^W (\nabla^F \tilde \rho)) (\tilde \varphi_{s(t)}(p)) .
\end{equation*}
Since $(\varphi^W_t)^{\ast} \nu = \e^{\tau t} \nu$, we obtain
\begin{align*}
\dv (\mathrm{d} \varphi_t^W (\nabla^F \tilde \rho) ) \nu &= \left. \frac{d}{ds} (\varphi^W_t \circ \tilde \varphi_s)^{\ast} \nu \right\vert_{s=0} = \left. \frac{d}{ds} \tilde \varphi_s^{\ast} (\varphi^W_t)^{\ast} \nu \right\vert_{s=0} \\ &= \e^{\tau t} \left. \frac{d}{ds} \tilde \varphi_s^{\ast} \nu \right\vert_{s=0} = \e^{\tau t} \dv (\nabla^F \tilde \rho) \nu .
\end{align*}
So
\begin{equation}
\label{LaplacianZermelo}
\triangle^Z \rho (\varphi_t (p)) = \dv W (\tilde \varphi_{s(t)} (p)) + s^{\prime\prime}(t) + s^{\prime}(t) \e^{\tau t} \triangle^F \tilde \rho (\tilde \varphi_{s(t)} (p)) .
\end{equation}

Using that $\dv W$ is constant again, we conclude $\triangle^Z \rho$ is constant on the fiber $\rho\vert_{\mathcal{U}}^{-1}(t) = \{\varphi_t (p) : p \in L\} \cap \mathcal{U}$ if only if $\triangle^F \tilde \rho$ is constant on the fiber $\tilde \rho\vert_{\tilde{\mathcal{U}}}^{-1} (s(t)) = \{\tilde \varphi_{s(t)} (p) : p \in L\} \cap  \tilde{\mathcal{U}}$. The remaining assertion is clear.
\end{proof}

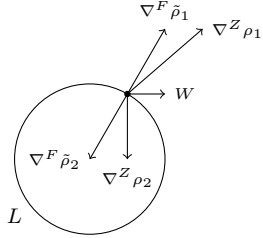
\begin{figure}[h]
\centering
\begin{tikzpicture}
\draw (0,0) circle [radius=1];
\node at (-1,-3/4) {\footnotesize $L$};
\filldraw (1/2,1.73205/2) circle (1pt) coordinate (p);
\path (p)  ++(60:1) coordinate (f1);
\path (p)  ++(60:-1) coordinate (f2);
\draw [->] (p) -- (f1) node[above] {\tiny $\nabla^F \tilde \rho_1$};
\draw[->] (p) -- (f2) node[below,left] {\tiny $\nabla^F \tilde \rho_2$};
\path (f1)  ++(0:1/2) coordinate (z1);
\path (f2)  ++(0:1/2) coordinate (z2);
\draw [->] (p) -- (z1) node[above,right] {\tiny $\nabla^Z \rho_1$};
\draw[->] (p) -- (z2) node[below] {\tiny $\nabla^Z \rho_2$};
\path (p) ++(0:1/2) coordinate (w);
\draw[->] (p) -- (w) node[right] {\tiny $W$};
\end{tikzpicture}
\caption{Choices for $\nabla^F \tilde \rho$ and $\nabla^Z \rho$.}
\label{gradfig}
\end{figure}

\begin{remark}
\label{gradhypothesis}
Let $\rho$ be a $Z$-distance function and suppose $L$ is a regular fiber of $\rho$. For each $p \in L$, $\nabla^Z \rho (p)$ is $Z$-orthonormal to $L$. By Lemma~\ref{Zermelofundtensor}, $\nabla^Z \rho (p) - W(p)$ is $F$-orthonormal to $L$. If $\tilde \rho$ is an $F$-distance function such that $L$ is a regular fiber of $\tilde \rho$, then $\nabla^F \tilde \rho\vert_L$ must be $F$-orthonormal to $L$. At each point $p \in L$, there are exactly two choices for $\nabla^F \tilde \rho (p)$. We pick $\tilde \rho$ so that $\nabla^F \tilde \rho (p) = \nabla^Z (p) - W(p)$, which translates to the hypothesis $\nabla^Z \rho \vert_L = (\nabla^F \tilde \rho + W)\vert_L$ in Theorem~\ref{Zermeloisop}. In particular, when $(M, Z)$ is a regular Randers space and $L$ is a closed hypersurface, $\nabla^Z \rho$ and $\nabla^F \tilde \rho$ lie on the same side of $L$; see Figure~\ref{gradfig}.

From a constructive perspective, we use the navigation process to deform an $F$-distance function $\tilde \rho$ into a $Z$-distance function $\rho$, as revealed in the proof of Theorem~\ref{Zermeloisop} by the equation $\varphi_t (p) = \varphi^W_t \circ \tilde \varphi_{s(t)} (p)$. In more details, let $\tilde \rho$ be an $F$-distance function and $L$ a regular fiber of $\tilde \rho$. If $\nabla^F \tilde \rho + W$ is $Z$-admissible on some tubular neighborhood of $L$, then $(\nabla^F \tilde \rho + W)\vert_L$ is $Z$-othonormal to $L$. For each $p \in L$ and $t$ on some neighborhood of zero (depending on $p$) the function $\rho$ defined by $\rho(\varphi^W_t \circ \tilde \varphi_{s(t)} (p)) = t$ is a $Z$-distance function by Theorem~\ref{Zermelogeod}. The hypersurface $L$ is a regular fiber of $\rho$ and $\nabla^Z \rho \vert_L = (\nabla^F \tilde \rho + W)\vert_L$. Conversely, Theorem~\ref{Zermelogeod} implies that any $Z$-distance function is locally obtained in this way.
\end{remark}

To conclude this section, we will show the hypothesis $\dv W$ constant is not an overreach. First, we prove that if $W$ is an $F$-homothetic vector field and $\nu$ is the Busemann-Hausdorff volume form associated with $Z$, then $\dv W$ is constant. But since the Busemann-Hausdorff volume form for $Z$ coincides with the one associated with $F$ by Proposition~\ref{ZermeloBH}, the following is enough.

\begin{lemma}
\label{divBH}
Let $(M, F, \nu_{BH})$ be a Finsler space with the Bussmann-Haudorff volume form. If $X$ is a homothetic vector field with $\varphi_t^{\ast} F = \e^{-\sigma t} F$, where $\varphi$ is the flow of $X$, then
\begin{equation*}
\varphi_t^{\ast} \nu_{BH} = \e^{-n \sigma t} \nu_{BH} ;
\end{equation*}
i.e. $\dv (W) = -n \sigma$.
\end{lemma}
\begin{proof}
Let $\{e_i\}_{i=1}^n$ be an oriented basis for $T_pM$, and let $\{\theta^i\}_{i=1}^n$ be the dual basis for $T_pM^{\ast}$. By definition, $(\nu_{BH})_p := \mu_{BH}(p) \theta^1 \wedge \cdots \wedge \theta^n$, with $\mu_{BH}$ given by Equation~\ref{BH}. Similarly, take an oriented basis $\{\tilde e_i\}_{i=1}^n$ for $T_{\varphi_t (p)}M$ and its dual basis $\{\tilde \theta^i\}_{i=1}^n$ for $T_{\varphi_t (p)}M^{\ast}$. Now, $(\varphi_t^{\ast} \nu_{BH})_p = \mu_{BH}(\varphi_t(p)) \varphi_t^{\ast} \tilde \theta^1 \wedge \cdots \wedge \varphi_t^{\ast} \tilde \theta^n$.

Since $\varphi_t$ is a local diffeomorphism, the differential map $\mathrm{d} \varphi_t : T_pM \to T_{\varphi_t(p)}M$ is invertible. Let $J_{\varphi_t}$ be the matrix representation of $\mathrm{d} \varphi_t$ with respect to the bases $\{e_i\}_{i=1}^n$ and $\{\tilde e_i\}_{i=1}^n$.

On one hand, $\varphi_t^{\ast} \tilde \theta^i = (J_{\varphi_t})^i_j \theta^j$. So $\varphi_t^{\ast} \tilde \theta^1 \wedge \cdots \wedge \varphi_t^{\ast} \tilde \theta^n = \det (J_{\varphi_t}) \theta^1 \wedge \cdots \wedge \theta^n$.

On the other hand, if $\tilde y^j := \e^{\sigma t} y^i (J_{\varphi_t})^j_i$, then
\begin{equation*}
F(\tilde y^j \tilde e_j) = F(\mathrm{d} \varphi_t (\e^{\sigma t} y^i e_i)) = \e^{-\sigma t} F(\e^{\sigma t} y^i e_i) = F(y^i e_i) ,
\end{equation*}
where the middle equality is the homothety condition for $X$ and the last follows from homogeneity. So $(\tilde y^j) \in B^n_{\varphi_t(p)}$ if and only if $(y^i) \in B^n_p$, and we have
\begin{align*}
\vl (B^n (\varphi_t(p))) &= \int_{B^n (\varphi_t(p))} \mathrm{d} \tilde y^1 \cdots \mathrm{d} \tilde y^n = \int_{B^n (p)} \e^{n \sigma t} \det(J_{\varphi_t}) \mathrm{d} y^1 \cdots \mathrm{d} y^n \\ &= \e^{n \sigma t} \det(J_{\varphi_t}) \vl(B^n (p))
\end{align*}
Hence, $\mu_{BH} (\varphi_t (p)) = \e^{-n \sigma t} \det(J_{\varphi_t})^{-1} \mu_{BH}(p)$. 

We conclude
\begin{equation*}
(\varphi_t^{\ast} \nu_{BH})_p = \e^{-n \sigma t} (\nu_{BH})_p .
\end{equation*}
\end{proof}

Next, if $W$ is $F$-homothetic and $\nu$ is the Holmes-Thompson volume form associated with $F$, then $\dv W$ is also constant.

\begin{lemma}
\label{divHT}
Let $(M, F, \nu_{HT})$ be a regular Finsler space with the Holmes-Thompson volume form. If $X$ is a homothetic vector field with $\varphi_t^{\ast} F = \e^{-\sigma t} F$, where $\varphi$ is the flow of $X$, then
\begin{equation*}
\varphi_t^{\ast} \nu_{HT} = \e^{-n \sigma t} \nu_{HT} ;
\end{equation*}
i.e. $\dv (W) = -n \sigma$.
\end{lemma}
\begin{proof}
Pick bases $\{e_i\}_{i=1}^n$, $\{\theta^i\}_{i=1}^n$, $\{\tilde e_i\}_{i=1}^n$, and $\{\tilde \theta^i\}_{i=1}^n$ as above.

By definition, $(\nu_{HT})_p := \mu_{HT}(p) \theta^1 \wedge \cdots \wedge \theta^n$, where $\mu_{HT}$ is given by Equation~\ref{HT}; similarly, $(\varphi_t^{\ast} \nu_{HT})_p = \mu_{HT}(\varphi_t(p)) \varphi_t^{\ast} \tilde \theta^1 \wedge \cdots \wedge \varphi_t^{\ast} \tilde \theta^n$.

Take $\tilde y^j := \e^{\sigma t} y^i (J_{\varphi_t})^j_i$. Once more, $(\tilde y^j) \in B^n (\varphi_t(p))$ if and only if $(y^i) \in B^n (p)$.

Let $\tilde y := \tilde y^j \tilde e_j$ and $y := y^i e_i$. Then $\tilde y = \mathrm{d} \varphi_t (\e^{\sigma t} y)$ and, by the homothety of $X$, $\varphi_t^{\ast} g = \e^{-2 \sigma t} g$. So
\begin{equation*}
g_{\tilde y} (\tilde e_i, \tilde e_j) = \e^{-2 \sigma t} g_{\e^{\sigma t} y} ((J_{\varphi_t}^{-1})^k_i e_k, (J_{\varphi_t}^{-1})^l_j e_l) = \e^{-2 \sigma t} g_{y} ((J_{\varphi_t}^{-1})^k_i e_k, (J_{\varphi_t}^{-1})^l_j e_l) ,
\end{equation*}
where the last equality follows from homogeneity of the fundamental tensor. Thus,
\begin{equation*}
g_{ij}(\tilde y) = \e^{-2 \sigma t} (J_{\varphi_t}^{-1})^k_i (J_{\varphi_t}^{-1})^l_j g_{kl}(y) ,
\end{equation*}
and we obtain $\det (g_{ij}(\tilde y)) = \e^{-2n \sigma t} \det (J_{\varphi_t})^{-2} \det (g_{ij}(y))$.

Hence,
\begin{equation*}
\int_{B^n (\varphi_t(p))} \det (g_{ij}(\tilde y)) \mathrm{d} \tilde y^1 \cdots \mathrm{d} \tilde y^n = \e^{-n \sigma t} \det (J_{\varphi_t})^{-1} \int_{B^n (p)} \det (g_{ij}(y)) \mathrm{d} y^1 \cdots \mathrm{d} y^n ;
\end{equation*}
i.e. $\mu_{HT}(\varphi_t(p)) = \e^{-n \sigma t} \det (J_{\varphi_t})^{-1} \mu_{HT}(p)$.

Therefore,
\begin{equation*}
(\varphi_t^{\ast} \nu_{HT})_p = \e^{-n \sigma t} (\nu_{HT})_p .
\end{equation*}
\end{proof}

\begin{remark}
\label{Finslervol}
Suppose $\nu^F$ is a volume form associated with $(M, F)$ that satisfies
\begin{enumerate}
\item $\nu^F$ is the Euclidean volume form  if $(M, F)$ is the Euclidean space;
\item $\int \varphi^{\ast} \nu^{\tilde F} = \int \nu^F$ if $\varphi: (M, F) \to (\tilde M, \tilde F)$ is an isometry;
\item $\int \nu^{F} \leq \int \nu^{\tilde F}$ if $F \leq \tilde F$.
\end{enumerate}
Then Theorem~\ref{Zermeloisop} applies when $W$ is an arbitrary $F$-homothetic vector field, for $\dv W$ will be constant. The conditions ensure the volume notion is Euclidean compatible and monotone with respect to the Finsler structure. In a Riemannian manifold, such volume form is uniquely determined, but there are many definitions in a general Finsler space, such as the Busemann-Hausforff and the Holmes-Thompson volume forms; we may call $\nu^F$ a \emph{Finslerian volume form}, e.g. \cite[5.5.3]{burago2001course}.
\end{remark}
\section{Proposition~\ref{Laplacediv}}

In order to prove Proposition~\ref{Laplacediv}, we establish an auxiliary lemma.

\begin{lemma}
\label{divgrad}
Let $(M,F,\nu)$ be any conic Finsler m space. If $f: \mathcal{U} \subset M \to \mathbb{R}$ is a smooth function, then at regular admissible points
\begin{equation*}
\dv \left( \frac{\nabla f}{F(\nabla f)} \right) = \frac{1}{F(\nabla f)} \left( \triangle f - \hs f \left( \frac{\nabla f}{F(\nabla f)}, \frac{\nabla f}{F(\nabla f)} \right) \right) .
\end{equation*}
\end{lemma}
\begin{proof}
By divergence properties,
\begin{equation*}
\dv \left( \frac{\nabla f}{F(\nabla f)} \right) = \frac{1}{F(\nabla f)} \dv (\nabla f) - \frac{1}{F^2(\nabla f)} \mathrm{d} ( F(\nabla f) ) (\nabla f).
\end{equation*}
On the other hand,
\begin{equation*}
F(\nabla f) \mathrm{d} ( F(\nabla f) ) (\nabla f) = \frac{1}{2} (\nabla f) ( F^2(\nabla f) ) = \hs f (\nabla f, \nabla f) ,
\end{equation*}
where the last equality follows from homogeneity of the fundamental and Cartan tensors together with almost metric compatibility of the Chern connection. Combine both equations to finish the proof.
\end{proof}

\begin{proof}[Proof of Proposition~\ref{Laplacediv}]
By Lemma~\ref{Zermelograd}(a) and divergence linearity,
\begin{equation*}
\dv \left( \frac{\nabla^Z f}{Z(\nabla^Z f)} \right) = \dv \left( \frac{\nabla^F f}{F(\nabla^F f)} \right) + \dv W.
\end{equation*}
Now Equation~\ref{Laplacefunc} follows from Lemma~\ref{divgrad}.
  
When $f$ is a $Z$-distance function, $Z( \nabla^Z f ) = 1$ and $\hs^Z f( \nabla^Z f, \nabla^Z f ) = 0$. So
\begin{equation*}
\triangle^Z f = \frac{1}{Z(\nabla^Z f)} \left( \triangle^Z f - \hs^Z f \left( \frac{\nabla^Z f}{Z(\nabla^Z f)}, \frac{\nabla^Z f}{Z(\nabla^Z f)} \right) \right) ,
\end{equation*}
and we conclude Equation~\ref{Laplacedistfunc} from Equation~\ref{Laplacefunc}.

If, additionally, $W$ is tangent to the fibers of $f$, then $\nabla^F f = \nabla^Z f - W$ by use of Lemma~\ref{Zermelofundtensor}. So $f$ is an $F$-distance function and, similarly to the above argument, we obtain Equation~\ref{Laplacedist}.
\end{proof}

\begin{corollary}
\label{Zermeloisopproj}
Let $(M, Z, \nu)$ be a Finsler m space with a conic Zermelo metric $Z$ of navigation data $(F,W)$. If $f : \mathcal{U} \subset M \to \mathbb{R}$ is a $Z$-admissible function, then the following are equivalent:
\begin{enumerate}
\item $W$ and $\dv W$ are $f$-projectable;
\item $f$ is $(F,\nu)$-isoparametric if and only if it is $(Z,\nu)$-isoparametric.
\end{enumerate}
\end{corollary}
\begin{proof}
Follows from Corollary~\ref{ZFtransnormal}, Lemma~\ref{Hessproj}, and Equation~\ref{Laplacefunc}.
\end{proof}

To some extent, Corollary~\ref{Zermeloisopproj} adapts Theorem~\ref{Zermeloisop} to \emph{one} arbitrary $Z$-admissible function $f$; in return, the conditions are more restrictive. Under the hypothesis of Theorem~\ref{Zermeloisop}, $W$ is an $F$-homothetic vector field such that $\dv W$ is constant, so $f$ is equivalently isoparametric for $(F, \nu)$ and $(Z, \nu)$ if and only if $W$ is $f$-projectable. Moreover, by Lemma~\ref{divBH}, if $W$ is \emph{any} $F$-homothetic vector field, then $f$ is simultaneously $(Z, \nu_{BH})$- and $(F, \nu_{BH})$-isoparametric if and only if $W$ is $f$-projectable; more generally, $\nu$ may be any Finslerian volume form associated with $F$ by Remark\ref{Finslervol}. Notably, this result extends \cite[Theorem~4.5]{xu2018isoparametric} for Randers spheres, and \cite[Theorem~3.9]{he2020classifications} for Randers spaces. In summary, the flow of $W$ limits which functions can be isoparametric with respect to both $(F, \nu)$ and $(Z, \nu)$; but conversely, $W$ may be restricted as a function $f$ is fixed in the construction of examples and counterexamples; cf. \cite[Example~3.4]{dong2022isoparametric}.

Consider, especially, a conic Randers manifold $(M, R)$ with navigation data $(h, W)$, where $(M, h)$ is a nonnegative Riemannian space form. Suppose $f$ is an $R$-transnormal function whose fibers are connected, and $W$ is an $f$-projectable vector field. By Corollary~\ref{ZFtransnormal}, $f$ is $h$-transnormal. But $f$ is $h$-transnormal if and only if it is $h$-isoparametric by \cite[Theorem~1.5(1)]{miyaoka2013transnormal} and the assumption of fiber connexity -- otherwise, the regular fibers of $f$ may form a family of $h$-isoparametric hypersurfaces without $f$ itself being an $h$-isoparametric function; see \cite[Example~5.1]{alexandrino1997hipersuperficies} or first footnote in \cite{thorbergsson1999survey}. Hence, $f$ is $(R, \nu_{BH})$-isoparametric if and only if $\dv W$ is $f$-projectable by Corollary~\ref{Zermeloisopproj}. When $W$ is $h$-homothetic, we use Lemma~\ref{divBH} to conclude $f$ \emph{is} $(R, \nu_{BH})$-isoparametric.
\section{Nonlinear Mean Curvature}

The nonlinear mean curvature of a hypersurface in a Finsler m space has been defined independently from different perspectives in at least two accounts, by Bellettini and Paolini \cite{bellettini1996anisotropic}, and by Shen \cite{shen1997curvature}. In the latter, the mean curvature is proposed constructively from the variation of the volume induced to a hypersurface by a unit normal vector field, and proved to coincide with the nonlinear Laplacian of a distance function, as per next result.

\begin{proposition}[Proposition~14.2, \cite{shen1997curvature}]
Let $(M, F, \nu)$ be a Finsler m space, and let $\rho$ be a smooth distance function on an open subset $\mathcal{U} \subset M$. The mean curvature $\Pi_{\xi}$ of a smooth level hypersurface $L := \rho^{-1}(c)$ with respect to $\xi := \nabla \rho$ is
\begin{equation*}
\Pi_{\xi (p)} = \triangle \rho (p) , \; \forall p \in L .
\end{equation*}
\end{proposition}

Bellettini and Paolini also demonstrated that their notion of mean curvature coincides (up to sign) with the Laplacian of a distance function; \cite[Lemma~3.2]{bellettini1996anisotropic}. In face of these results, we choose to extend the mean curvature to conic Finsler spaces by the Laplacian of a distance function, and then establish the variation of the volume, as in \cite[Theorem~5.1]{bellettini1996anisotropic} and \cite[Section~14]{shen1997curvature}.

Let $(M, F, \nu)$ be a conic Finsler m space, and $L \subset M$ be an admissible hypersurface. For each $p \in L$, there is a distance function $\rho$ defined on an open subset $\mathcal{U} \subset M$ around $p$ such that $L \cap \mathcal{U}$ is a regular level set of $\rho$ and $\nabla \rho = \xi$ along $L \cap \mathcal{U}$. The \emph{(nonlinear) mean curvature} -- or \emph{$(F,\nu)$-mean curvature} -- of $L$ at $p$ in the direction of $\xi$ is
\begin{equation*}
\Pi_{\xi(p)} := \triangle \rho (p) .
\end{equation*}
The existence of $\rho$ and the independence of $\Pi_{\xi}$ on the choice of $\rho$ for conic Finsler spaces are both consequences of Lemma~\ref{existofdistfunc}.

\begin{lemma}
\label{meancurveq}
Let $(M, F, \nu)$ be a conic Finsler m space, and let $f$ be an admissible function on an open subset $\mathcal{U} \subset M$. The mean curvature of a regular level hypersurface $L := f^{-1}(c)$ with respect to $\frac{\nabla f}{F(\nabla f)}$ is
\begin{equation*}
\Pi_{\frac{\nabla f(p)}{F(\nabla f(p))}} = \frac{1}{F(\nabla f(p))}\left( \triangle f(p) - \hs f \left(\frac{\nabla f(p)}{F(\nabla f(p))}, \frac{\nabla f(p)}{F(\nabla f(p))}\right) \right) , \; \forall p \in L .
\end{equation*}
\end{lemma}
\begin{proof}
If $L:= f^{-1}(c)$ is a regular level hypersurface, then $\frac{\nabla f}{F(\nabla f)}$ is an admissible orthonormal vector field along $L$. By Lemma~\ref{existofdistfunc}, there is a locally defined distance function $\rho$ around a fixed point $p \in L$ such that $\nabla \rho\vert_{\mathcal{U} \cap L} = \frac{\nabla f}{F(\nabla f)}\vert_{\mathcal{U} \cap L}$. So we conclude the proof by Lemmas~\ref{divhypersurface} and~\ref{divgrad}.
\end{proof}

Clearly, for Riemannian manifolds, the nonlinear mean curvature associated with any Finslerian volume form -- such as the Busemann-Hausdorff or the Holmes-Thompson -- coincides with the Riemannian mean curvature.

\begin{proposition}
Let $(M, F, \nu)$ be a conic Finsler m space, and let $L \subset M$ be an admissible hypersurface for which $\xi$ is an admissible unit normal vector field. For each $p \in L$, let $\gamma_{\xi(p)}(t)$ denote the unit speed geodesic with $\gamma_{\xi(p)}(0) = p$ and $\gamma^{\prime}_{\xi(p)}(0) = \xi (p)$. Suppose $\varphi: L \to \mathbb{R}$ is a positive function with compact support, and $\varepsilon > 0$ is small enough for the function $f_t (p) := \gamma_{\xi(p)} (\varphi(p) t)$ on $\supp \varphi$ to be an immersion for all $t \in (-\varepsilon, \varepsilon)$. If the volume form in $L_t := \im f_t$ induced by $X( f_t(p) ) := \gamma^{\prime}_{\xi(p)} (\varphi(p) t)$ is $\nu_t := i_X \nu$, and $\vl (L_t) :=\int_{L_t} \nu_t$, then
\begin{equation*}
\left. \frac{d}{dt} \vl (L_t) \right\vert_{t=0}	= \int_{L_0} \varphi(p) \Pi_{\xi(p)} \nu_0 .
\end{equation*}
\end{proposition}
\begin{proof}
Let $\{ E_1, \ldots, E_n \}$ be a local frame for $TM$ with $E_1 := X$. Consider the function $\Phi (\gamma_{\xi(p)} (t)) := \varphi(p)$. For all $p \in L_0$,
\begin{equation*}
\mathrm{d} \Phi_p (X) = \frac{1}{\varphi(p)} \left. \frac{d}{dt} \Phi (f_t(p)) \right\vert_{t=0} = \frac{1}{\varphi(p)} \left. \frac{d}{dt} \varphi(p) \right\vert_{t=0} = 0 .
\end{equation*}
Moreover, $(\mathcal{L}_X \nu)_p = \Pi_{\xi(p)} \nu_p$ by use of Lemmas~\ref{existofdistfunc} and ~\ref{divhypersurface}. So,
\begin{equation*}
(\mathcal{L}_{\Phi X} \nu)_p (E_1, \ldots, E_n) = \varphi(p) \Pi_{\xi(p)} \nu_p (E_1, \ldots, E_n) .
\end{equation*}
Since $\{E_2(p), \ldots, E_n(p)\}$ is a basis for $T_p L$,
\begin{align*}
\left. \frac{d}{dt} (f_t^{\ast} \nu_t) \right\vert_{t=0} (E_2, \ldots, E_n) &= (\mathcal{L}_{\Phi X} \nu_t) (E_2, \ldots, E_n) = (\mathcal{L}_{\Phi X} \nu) (E_1, \ldots, E_n) \\
&= \varphi(p) \Pi_{\xi(p)} \nu (E_1, \ldots, E_n) = \varphi(p) \Pi_{\xi(p)} \nu_0 (E_2, \ldots, E_n) .
\end{align*}

Therefore,
\begin{equation*}
\left. \frac{d}{dt} \vl (L_t) \right\vert_{t=0}	= \int_{L_0} \left. \frac{d}{dt}(f_t^{\ast} \nu_t) \right\vert_{t=0} = \int_{L_0} \varphi(p) \Pi_{\xi(p)} \nu_0 .
\end{equation*}
\end{proof}

For any Finsler m space $(M, F, \nu)$, the regular fibers of a transnormal function $f$ form a family of isoparamentric hypersurfaces if and only if they have constant nonlinear mean curvature with respect to the unitary normal vector field $\frac{\nabla f}{F(\nabla f)}$, by \cite[Theorem~1.1]{he2016isoparametric}. We may rewrite and proof this fact for conic Finsler spaces in the following way.

\begin{theorem}
\label{isopmeancurv}
Let $(M, F, \nu)$ be a conic Finsler m space. Suppose $f: \mathcal{U} \subset M \to \mathbb{R}$ is a transnormal function. Then $f$ is $(F, \nu)$-isoparametric if only if the function $p \mapsto \Pi_{\frac{\nabla f (p)}{F(\nabla f (p))}}$ is $f$-projectable on $\{p \in \mathcal{U} : \mathrm{d} f(p) \neq 0\}$.
\end{theorem}
\begin{proof}
By Lemma~\ref{meancurveq}, for all regular level hypersurfaces we have
\begin{equation*}
\Pi_{\frac{\nabla f}{F(\nabla f)}} = \frac{1}{F(\nabla f)}\left( \triangle f - \hs f \left(\frac{\nabla f}{F(\nabla f)}, \frac{\nabla f}{F(\nabla f)}\right) \right).
\end{equation*}

If $f$ is transnormal, then $F(\nabla f)$ and $\hs f\left( \frac{\nabla f}{F(\nabla f)}, \frac{\nabla f}{F(\nabla f)} \right)$ are $f$-projectable by definition and Lemma~\ref{Hessproj}. Therefore, $\triangle f$ is $f$-projectable if and only if $\Pi_{\frac{\nabla f}{F(\nabla f)}}$ is $f$-projectable.
\end{proof}

Next, we prove Corollary~\ref{Zermelomean}, which relates the nonlinear mean curvature with respect to the Busemann-Hausdorff volume form for a conic Zermelo space $(M,Z)$ and the Finsler manifold $(M, F)$ from the navigation data.

\begin{proof}[Proof of Corollary~\ref{Zermelomean}]
If $\xi$ is an admissible unit normal vector field with respect to $Z$ along a hypersurface $L \subset M$, then $\xi - W$ is a unit normal vector field with respect to $F$ along $L$ by Lemma~\ref{Zermelofundtensor}.

For each $p \in L$, take a $Z$-distance function $\rho$ on some open neighborhood $\mathcal{U} \subset M$ of $p$ such that $L \cap \mathcal{U}$ is a regular level set of $\rho$ and $\nabla^Z \rho (p) = \xi (p)$. By Lemma~\ref{Zermelograd}(a), at all regular admissible points of $\rho$ we have
\begin{equation*}
\frac{\nabla^F \rho}{F (\nabla^F \rho)} = \nabla^Z \rho - W .
\end{equation*}
So the $(F, \nu^F_{BH})$-mean curvature of $L$ in the direction of $\nabla^Z \rho - W$ is, by Lemma~\ref{meancurveq},
\begin{equation*}
\Pi^F_{\nabla^Z \rho - W} = \frac{1}{F (\nabla^F \rho)} (\triangle^F \rho - \hs^F \rho(\nabla^Z \rho - W, \nabla^Z \rho - W) ) .
\end{equation*}

By Proposition~\ref{ZermeloBH}, the Busemann-Hausdorff volume form for $Z$ coincides with the one associated with $F$. Therefore, by use of Equation~\ref{Laplacedistfunc}, the $(Z,\nu_{BH})$-mean curvature of $L$ with respect to $\nabla^Z \rho$ is
\begin{equation*}
\Pi^Z_{\nabla^Z \rho} = \Pi^F_{\nabla^Z \rho - W} + \dv W .
\end{equation*}
Evaluate the expression at $p$ to conclude the proof.
\end{proof}

When the conic Zermelo metric arrises from homothetic navigation, $\dv W$ is constant by Lemma~\ref{divBH}. So the nonlinear mean curvature with respect to $(Z, \nu_{BH})$ of an oriented hypersurface $L$ is constant if and only if the $(F, \nu_{BH})$-mean curvature of $L$ is constant. This result was first proved in \cite[Theorem~1.2]{qian2022hypersurfaces}. It is also a consequence of Corollary~\ref{Zermeloisopproj} in light of Theorem~\ref{isopmeancurv}. In fact, an alternative proof follows from Equation~\ref{LaplacianZermelo}, which holds for any conic Zermelo pseudo-metric.

We finish this section with an explicit use of Corollary~\ref{Zermelomean}.

\begin{example}
\label{nonlinearmeanex}
Let $(\mathbb{R}^n, h)$ be the Euclidean space for $n \geq 2$ and let $W$ be the smooth vector field $W(x) = - x$, $\forall x \in \mathbb{R}^n$. Consider the conic Randers metric $R$ with navigation data $(h, W)$. The restriction of $R$ to the open unit ball $\mathbb{B}^n := \{ x \in \mathbb{R}^n : \vert x \vert < 1 \}$ coincides with the Funk metric.

Define the function $f(x) := \vert x \vert$, $\forall x \in \mathbb{B}^n \setminus 0$, and let $S_r := f^{-1} (r)$, with $0 < r < 1$. The vectors $\pm \nabla^h f (x)$ are the orthonormal vectors to $S_r$ at $x$ with respect to $h$. Moreover, since $f$ is an $h$-distance function, the Riemannian mean curvature of $S_r$ in the direction of $\pm \nabla^h f$ is
\begin{equation*}
\Pi^h_{\pm \nabla^h f} = \pm \triangle^h f = \pm \frac{(n-1)}{r} .
\end{equation*}
Hence, $f$ is an $h$-isoparametric function. Once $W$ and $\dv W$ are $f$-projectable, Corollary~\ref{Zermeloisopproj} implies $f$ is $(R, \nu_{BH})$-isoparametric.

By Lemma~\ref{Zermelograd}(a), the vectors $\pm \nabla^h f (x) + W (x)$ are the orthonormal vectors to $S_r$ at $x$ with respect to $R$. So by Corolary~\ref{Zermelomean}, the nonlinear mean curvature for $(R, \nu_{BH})$ of $S_r$ in the direction of $\pm \nabla^h f + W$ is
\begin{equation*}
\Pi^R_{\pm \nabla^h f + W} = \Pi^h_{\pm \nabla^h f} + \dv W = \pm \frac{(n-1)}{r} - n .
\end{equation*}
In particular, the $(R, \nu_{BH})$-mean curvature of $S_{1-\frac{1}{n}}$ in the direction of $\nabla^h f + W$ is null!
\end{example}
\section{Linear Mean Curvature}

First, we recall the definition of linear mean curvature presented by Shen in \cite{shen1998finsler} with minor adjustments to our writing choices.

Let $(M^n, F)$ be a Finsler manifold, and let the inclusion map $i: L^m \to (M^n, F)$ be an immersion. Clearly, a smooth vector field on $(L^m, i^{\ast} F)$ is identified as a smooth vector field along $L$ in $(M^n, F)$ by the bundle map $i_{\ast}: TL \to TM$. Conversely, there is a canonical projection $P_{i_{\ast}}: TM \to TL$ so that each smooth vector field $X$ on $M$ corresponds to a smooth vector field $P_{i_{\ast}} (X)$ on $L$; refer to \cite[Section~5]{shen1998finsler} for the explicit definition.

Suppose $(L^m, i^{\ast} F)$ is oriented and consider a precompact open subset $\mathcal{U} \subset L$. Take a smooth family of immersions $f_t: L^m \to (M^n, F)$ such that $f_0 = i$, and $f_t \equiv i$ outside the closure of $\mathcal{U}$. The metric induced by $f_t$ on $L$ is $F_t := f_t^{\ast} F$; in particular, $F_0 = i^{\ast} F$. For each $p \in L$, choose an arbitrary oriented basis $\{ e_i \}_{i=1}^{m}$ for $T_p L$, and its dual basis $\{\theta^i \}_{i=1}^{m}$ for $T^{\ast}_p L$. The Busemann-Hausdorff volume form for the induced metric $F_t$ is $\nu_{BH}^{F_t} = \mu_t \theta^1 \wedge \ldots \wedge \theta^m$, where $\mu_t$ is a simplified notation for the function $\mu_{BH}^{F_t}$ defined according to Equation~\ref{BH}.

Denote the variational vector field along $L$ by $X := \left. \frac{\partial}{\partial t} f_t \right\vert_{t=0}$. The \emph{linear mean curvature} of $L$ in the direction $X$ is
\begin{equation}
\label{linearmeancurv}
\mathcal{H} (X) := \left. \frac{d}{d t} \ln (\mu_t) \right\vert_{t=0} - \dv [P_{f_{\ast}} (X)] ,
\end{equation}
where $\dv$ is the divergence in $(L^m, \nu_{BH}^{F_0})$. This quantity arises from the variation by isometric immersions of the Busemann-Hausdorff volume form associated with the induced metrics on $L$.

\begin{theorem}[Theorem~1.2, \cite{shen1998finsler}]
\label{linearmeanvariation}
Let $L^m$ be an immersed submanifold of the Finsler space $(M^n, F)$ by the inclusion map $i: L \to M$. Suppose $\mathcal{U} \subset L$ is a precompact open subset, and $f_t: (L^m, F_t) \to (M^n, F)$ is a smooth variation of isometric immersions with $f_0 = i$ and $f_t \equiv i$ in $L \setminus \overline{\mathcal{U}}$. If the Busemann-Hausdorff volume form associated with $F_t$ is $\nu_{BH}^{F_t}$ and $\vl_t (\mathcal{U}) := \int_{\mathcal{U}} \nu_{BH}^{F_t}$, then
\begin{equation*}
\left. \frac{d}{d t} \vl_t (\mathcal{U}) \right\vert_{t=0} =  \int_{\mathcal{U}} \mathcal{H} (X) \nu_{BH}^{F_0} ,
\end{equation*}
where $X := \left. \frac{\partial}{\partial t} f_t \right\vert_{t=0}$ is the variational vector field.
\end{theorem}
 
We can now proceed to our results. Let $(M, R)$ be a Randers manifold with navigation data $(h, W)$ and suppose $L^m \subset M^n$ is an immersed submanifold.

\begin{proof}[Proof of Theorem~\ref{Randerslinearmean}]
Let $i: L^m \to (M^n, R)$ denote the inclusion map. For each precompact open subset $\mathcal{U} \subset L$, suppose $f_t: (L^m, R_t) \to (M^n, R)$ is a smooth variation of isometric immersions such that $f_0 = i$ and $f_t \equiv i$ in $L \setminus \overline{\mathcal{U}}$.

For $t$ small enough, $L_t := \im f_t \subset M$ is an immersed submanifold. Let $W^{\perp}(f_t(p))$ be the component of the vector field $W$ which is $h$-orthogonal to $L_t$ at $f_t(p)$. Define
\begin{equation*}
\lambda_t (p) := \frac{1}{1 - h(W^{\perp}(f_t(p)), W^{\perp}(f_t(p)))} \, , \forall p \in L \, .
\end{equation*}
Once $f_t$ gives an isometry between $(L, R_t)$ and $(L_t, R\vert_{TL_t})$, it follows from Corollary~\ref{RandersimmersionBH} that $\nu^{R_t}_{BH} = \lambda_t^{\frac{m}{2}} \nu^{h_t}$, where $h_t := f^{\ast}_t (h\vert_{TL_t \times TL_t})$. Denote the variational vector field by $X := \left. \frac{\partial}{\partial t} f_t \right\vert_{t=0}$. Now, by Theorem~\ref{linearmeanvariation}, we have that
\begin{align*}
\int_{\mathcal{U}} \mathcal{H}(X) \nu^{R_0}_{BH} &= \int_{\mathcal{U}} \left. \frac{d}{d t} \nu^{R_t}_{BH} \right\vert_{t=0} = \int_{\mathcal{U}} \left( \frac{m}{2} \lambda_0^{\frac{m}{2}-1} \left. \frac{d}{d t} \lambda_t \right\vert_{t=0} \nu^{h_0} + \lambda_0^{\frac{m}{2}} \left. \frac{d}{d t} \nu^{h_t} \right\vert_{t=0} \right) \\
&= \int_{\mathcal{U}} \left( \frac{m}{2} \lambda_0^{-1} \left. \frac{d}{d t} \lambda_t \right\vert_{t=0} + \Pi^h_X \right) \nu_{BH}^{R_0}
\end{align*}
Thus,
\begin{equation}
\label{EqA}
\mathcal{H}(X) = \Pi^{h}_{X} + \left. \frac{d}{d t} \ln \lambda_t^{\frac{m}{2}} \right\vert_{t=0} .
\end{equation}

Suppose $\mathbf{n}$ is an orthonormal vector field along $L$ with respect to $h$. For each $p \in L$, let $\beta_p(t)$ be the unit speed geodesic with respect to $h$ satisfying $\beta_p(0) = p$ and $\beta_p^{\prime}(0) = \mathbf{n}(p)$. Consider $f_t(p) := \beta_p(t)$ for $p \in \mathcal{U} \subset L$ and $t$ small enough. Then, $W^{\perp}(f_t(p)) = h(W(\beta_p(t)), \beta_p^{\prime}(t)) \beta_p^{\prime}(t)$, and so $\lambda_t(p) = \frac{1}{1-h(W(\beta_p(t)), \beta_p^{\prime}(t))^2}$. Hence,
\begin{equation}
\label{EqB}
\left. \frac{d}{d t} \ln \lambda_t^{\frac{m}{2}} \right\vert_{t=0} = \frac{m h(W, \mathbf{n}) h(D_{\mathbf{n}}W, \mathbf{n})}{1 - h(W, \mathbf{n})^2} =: \mathcal{B}(\mathbf{n}) ,
\end{equation}
where $D$ is the Levi-Civita connection for $h$ and we use $D_{\beta_p^{\prime}(t)} \beta_p^{\prime}(t) = 0$.

Put together Equations \ref{EqA} and \ref{EqB} to conclude the proof.
\end{proof}

\begin{example}
\label{linearmeanex}
Let $(M^n, h)$ be a Riemannian manifold with $n \geq 2$ and let $L \subset M$ be a hypersurface. Suppose $\mathbf{n}$ is an $h$-orthonormal vector field along $L$ and denote by $\beta_p(t)$ the unit speed geodesic with respect to $h$ satisfying $\beta_p(0) = p$ and $\beta_p^{\prime}(0) = \mathbf{n}(p)$. Assume there exists $\varepsilon > 0$ such that the map $(p, t) \mapsto  \beta_p(t)$ is a diffeomorphism of $L \times (-\varepsilon, \varepsilon)$ onto its image; e.g. when $L$ is compact.

Define the tubular neighborhoods $\mathcal{U}_{\delta} (L) := \{ \beta_p(t) : p \in L , \vert t \vert < \delta \} \supset L$ for any $0 < \delta \leq \varepsilon$. For each $p \in L$, choose a constant $\pi(p) > 0$ such that the function $c(\beta_p(t)) := \frac{2}{(n-1)}\sqrt{1 - \e^{t\Pi^h_{\mathbf{n}(p)}-\pi(p)}}$ is well-defined, smooth and positive on $\mathcal{U}_{\varepsilon} (L)$. By Tietze extension theorem, there is a smooth function $C$ on $M$ such that $C\vert_{\overline{\mathcal{U}_{\delta}(L)}} = c$ and $\im C \subset (0, 1)$ for some $0 < \delta < \varepsilon$.

Consider an extension $X$ on $M$ of the smooth vector field $\beta_p^{\prime}(t)$ on $\mathcal{U}_{\varepsilon}(L)$ such that $h(X, X) \leq 1$. For the smooth vector field $W:= C X$ we have that $h(W, W) < 1$. So the Finsler metric $R$ with navigation data $(h, W)$ is a regular Randers metric.

Since $\mathrm{d}c (\mathbf{n}) = -\frac{(1-c^2) \Pi^h_{\mathbf{n}}}{(n-1) c}$, it follows that
\begin{equation*}
\mathcal{B}(\mathbf{n}) = \frac{(n-1) h(c\mathbf{n}, \mathbf{n}) h(D_{\mathbf{n}}CX, \mathbf{n})}{1 - h(c\mathbf{n}, \mathbf{n})^2} = \frac{(n-1) c}{1-c^2} \mathrm{d}c(\mathbf{n}) = - \Pi^h_{\mathbf{n}} ,
\end{equation*}
where $D$ is the Levi-Civita connection for $h$ and $D_{\beta_p^{\prime}(t)} \beta_p^{\prime}(t) = 0$. Hence, $L$ is a BH-minimal surface of $(M, R)$ by Theorem~\ref{Randerslinearmean}.
\end{example}

In summary, any hypersurface $L \subset M$ which admits a tubular neighborhood and an orthonormal vector field with respect to some given Riemannian metric $h$ is a BH-minimal surface with respect to some Randers metric $R$ on $M$.

\begin{corollary}
\label{Randersisopproj}
Let $(M, R, \nu_{BH})$ be a Randers m space with navigation data $(h, W)$ and the Busemann-Hausdorff volume form. If a function $f$ is $h$- or $R$-transnormal, then the following are equivalent:
\begin{enumerate}
\item $W$ and $\dv(W)$ are $f$-projectable;
\item $f$ is $(R, \nu_{BH})$-isoparametric if and only if $p \mapsto \mathcal{H} (\nabla^h f(p))$ is $f$-projectable.
\end{enumerate} 
In case either is true, the singular fibers of $f$ are BH-minimal submanifolds when the function is $(R, \nu_{BH})$-isoparametric.
\end{corollary}
\begin{proof}
By Lemma~\ref{Zermelograd}(a), $\frac{\nabla^R f}{R(\nabla^R f)} = \frac{\nabla^h f}{\vert \nabla^h f \vert} + W$, where $\vert v \vert := \sqrt{h(v, v)}$, for all $v \in TM$. So the linear and the nonlinear mean curvatures of any regular level hypersurface of $f$ in $(M, R, \nu_{BH})$ are related by Corollary~\ref{Zermelomean} and Theorem~\ref{Randerslinearmean} acccording to the equation
\begin{equation}
\label{meancurvs}
\mathcal{H}\left( \frac{\nabla^h f}{\vert \nabla^h f \vert} \right) = \Pi^R_{\frac{\nabla^h f}{\vert \nabla^h f \vert} + W} - \dv W + \mathcal{B}\left( \frac{\nabla^h f}{\vert \nabla^h f \vert} \right) .
\end{equation}

By Lemma~\ref{Zermelograd}(b) and the hypothesis, $W$ is $f$-projectable if and only if $f$ is both $h$- and $R$-trasnormal. Given such equivalence, if either holds, then $\mathcal{B}\left( \frac{\nabla^h f}{\vert \nabla^h f \vert} \right)$ is $f$-projectable, because $W$ is $f$-projectable too. Moreover, Theorem~\ref{isopmeancurv} implies that $f$ is $(R, \nu_{BH})$-isoparametric if and only if the $(R, \nu_{BH})$-mean curvature in the direction of $\frac{\nabla^R f}{R(\nabla^R f)}$ is $f$-projectable over the regular fibers. By Equation~\ref{meancurvs} and the linearity of $\mathcal{H}$, we conclude $f$ is $(R, \nu_{BH})$-isoparametric if and only if $\mathcal{H}(\nabla^h f)$ is $f$-projectable exactly when $\dv W$ is also $f$-projectable.

The last affirmation is now proved similarly to the Riemannian case.
\end{proof}

Essentially, Corollary~\ref{Randersisopproj} gives a characterization for $(R, \nu_{BH})$-isoparametric functions with respect to the linear mean curvature of their regular fibers in some special Randers space. A naive question is whether the linear mean curvature could characterize isoparametric functions in any Finslerian manifold. But Corollary~\ref{Randersisopproj} also provides the negative answer even for Randers manifolds, according to the example below.

\begin{example}
\label{linearmeancounterex}
Let $(\mathbb{H}^n, h)$ be the hyperbolic space with $n \geq 2$ and let $\rho$ be a regular $h$-distance function on $\mathbb{H}^n$ that is \emph{not} $h$-isoparametric; refer to the proof of Theorem~1.5(3) in \cite{miyaoka2013transnormal} for the construction of such function. 

Choose a nonvanishing smooth vector field $X$ on $\mathbb{H}^n$ such that $X$ is tangent to the fibers of $\rho$ and $h(X, X) < 1$; its existence is ensured by the assumption $n \geq 2$ and the regularity of $\rho$. For the linear differential equation $X(c) + (\dv X)c + \triangle^h \rho = 0$, there exists a nontrivial solution $c: \mathcal{U} \subset \mathbb{H}^n \to \mathbb{R}$ such that $\vert c \vert < 1$, taking $\mathcal{U}$ small enough. So the smooth vector field $W := c X$ on $\mathcal{U}$ satisfies $h(W,W) < 1$ and
\begin{equation*}
\dv W = X(c) + (\dv X)c = -\triangle^h \rho .
\end{equation*}

Consider the regular Randers space $(\mathcal{U}, R)$ with navigation data $(h, W)$. Recall that the Busemann-Hausdorff volume form $\nu_{BH}$ associated with $R$ coincides with the Riemannian volume form $\nu^h$ on $\mathcal{U}$; vide Proposition~\ref{ZermeloBH}. Thus, the divergence for $(\mathcal{U}, R)$ is the same as the divergence with respect to $h$, denoted $\dv$ as above.

Since $W$ is tangent to the fibers of $\rho\vert_{\mathcal{U}}$, we have that $\mathbf{d} \rho (W) = 0$. So $\rho$ is an $R$-distance function by Lemma~\ref{Zermelograd}(b). Now, Equation~\ref{Laplacedist} implies
\begin{equation*}
\triangle^R \rho = \triangle^h \rho + \dv W = 0.
\end{equation*}
Hence, $\rho$ is $(R, \nu_{BH})$-isoparametric. However, by Theorem~\ref{Randerslinearmean}, $\mathcal{H}(\nabla^h\rho) = \triangle^h \rho$ is \emph{not} $\rho$-projectable.
\end{example}

The example suggests that the linear mean curvature is not enough to study isoparametric hypersurfaces; hence, it is not an appropriate notion to extend the definition of isoparametric hypersurfaces to general isoparametric submanifolds. Is it possible to find another notion to describe isoparametric submanifolds?

\bibliography{references}

\providecommand{\bysame}{\leavevmode\hbox to3em{\hrulefill}\thinspace}
\providecommand{\MR}{\relax\ifhmode\unskip\space\fi MR }
\providecommand{\MRhref}[2]{%
  \href{http://www.ams.org/mathscinet-getitem?mr=#1}{#2}
}
\providecommand{\href}[2]{#2}
\begin{thebibliography}{10}

\bibitem{alexandrino1997hipersuperficies}
Marcos~M Alexandrino, \emph{{Hipersuperf\'icies de n\'ivel de uma fun\c c\~ao
  transnormal}}, Master's thesis, Pontif\'icia Universidade Cat\'olica do Rio
  de Janeiro, 1997.

\bibitem{alexandrino2019finsler}
Marcos~M Alexandrino, Benigno~O Alves, and Hengameh~R Dehkordi, \emph{{On
  Finsler transnormal functions}}, Differential Geometry and its Applications
  \textbf{65} (2019), 93--107.

\bibitem{anastasiei1993absolute}
Mihai Anastasiei and Hiroaki Kawaguchi, \emph{{Absolute energy of a Finsler
  space can't be harmonique}}, Tensor, New Series \textbf{53} (1993), 108--114.

\bibitem{antonelli1993stochastic}
Peter~L Antonelli and Tomasz~J Zastawniak, \emph{{Stochastic calculus on
  Finsler manifolds and an application in biology}}, Nonlinear World \textbf{1}
  (1993), 149--171.

\bibitem{balan2007bh}
Vladimir Balan, \emph{{BH-mean curvature in Randers spaces with anisotropically
  scaled metric}}, {Proceedings of The International Conference ``Differential
  Geometry and Dynamical Systems'' (DGDS-2007)} (Constantin Udriste and
  Vladimir Balan, eds.), Geometry Balkan Press, 2008, pp.~34--39.

\bibitem{bao1996hodge}
David Bao and Brad Lackey, \emph{{A Hodge decomposition theorem for Finsler
  spaces}}, Comptes rendus de l'Acad{\'e}mie des sciences. S{\'e}rie 1,
  Math{\'e}matique \textbf{323} (1996), no.~1, 51--56.

\bibitem{bao1996eigenforms}
\bysame, \emph{{Special eigenforms on the sphere bundle of a Finsler
  manifold}}, Contemporary Mathematics \textbf{196} (1996), 67--78.

\bibitem{bao2004zermelo}
David Bao, Colleen Robles, and Zhongmin Shen, \emph{{Zermelo navigation on
  Riemannian manifolds}}, Journal of Differential Geometry \textbf{66} (2004),
  no.~3, 377--435.

\bibitem{bellettini1996anisotropic}
Giovanni Bellettini and Maurizio Paolini, \emph{{Anisotropic motion by mean
  curvature in the context of Finsler geometry}}, Hokkaido Mathematical Journal
  \textbf{25} (1996), no.~3, 537--566.

\bibitem{burago2001course}
Dmitri Burago, Yuri Burago, and Sergei Ivanov, \emph{{A Course in Metric
  Geometry}}, Graduate Studies in Mathematics, vol.~33, American Mathematical
  Society, 2001.

\bibitem{caponio2014wind}
Erasmo Caponio, Miguel~A Javaloyes, and Miguel S{\'a}nchez, \emph{{Wind
  Finslerian structures: from Zermelo's navigation to the causality of
  spacetimes}},
  \href{https://arxiv.org/abs/1407.5494}{\texttt{arXiv.1407.5494}}, 2014.

\bibitem{cartan1938familles}
{\'E}lie Cartan, \emph{{Familles de surfaces isoparam\'etriques dans les
  espaces \`a courbure constante}}, Annali di Matematica Pura ed Applicata
  \textbf{17} (1938), no.~1, 177--191.

\bibitem{cartan1939sur}
\bysame, \emph{{Sur des familles remarquables d'hypersurfaces
  isoparam\'etriques dans les espaces sph\'eriques}}, Mathematische Zeitschrift
  \textbf{45} (1939), no.~1, 335--367.

\bibitem{chakerian1996integral}
Gulbank~D Chakerian, \emph{{Integral geometry in Minkowski spaces}},
  Contemporary Mathematics \textbf{196} (1996), 43--50.

\bibitem{chen2022transnormal}
Yali Chen and Qun He, \emph{{Transnormal functions and focal varieties on
  Finsler manifolds}},
  \href{https://arxiv.org/abs/2203.09796}{\texttt{arXiv.2203.09796}}, 2022.

\bibitem{chen2023isoparametric}
\bysame, \emph{{The isoparametric functions on a class of Finsler spheres}},
  Differential Geometry and its Applications \textbf{86} (2023), 101970.

\bibitem{chi2020isoparametric}
Quo-Shin Chi, \emph{{Isoparametric hypersurfaces with four principal
  curvatures, IV}}, Journal of Differential Geometry \textbf{115} (2020),
  no.~2, 225--301.

\bibitem{cricsan2020finsler}
Adina~V Cri{\c{s}}an and Ion~V Vancea, \emph{{Finsler geometries from
  topological electromagnetism}}, The European Physical Journal C \textbf{80}
  (2020), no.~6, 1--12.

\bibitem{cui2014minimal}
Ningwei Cui, \emph{{On minimal surfaces in a class of Finsler $3$-spheres}},
  Geometriae Dedicata \textbf{168} (2014), no.~1, 87--100.

\bibitem{cui2017nontrivial}
Ningwei Cui and Yi-Bing Shen, \emph{{Nontrivial minimal surfaces in a
  hyperbolic Randers space}}, Mathematische Nachrichten \textbf{290} (2017),
  no.~4, 570--582.

\bibitem{cvetivc2012graphene}
Mirjam Cveti{\v{c}} and Gary~W Gibbons, \emph{{Graphene and the Zermelo optical
  metric of the BTZ black hole}}, Annals of Physics \textbf{327} (2012),
  no.~11, 2617--2626.

\bibitem{dasilva2011minimal}
Rosangela~Maria da~Silva and Keti Tenenblat, \emph{{Minimal surfaces in a
  cylindrical region of $\mathbb{R}^3$ with a Randers metric}}, Houston Journal
  of Mathematics \textbf{37} (2011), no.~3, 745--771.

\bibitem{dasilva2014helicoidal}
Ros{\^a}ngela~Maria da~Silva and Keti Tenenblat, \emph{{Helicoidal minimal
  surfaces in a Finsler space of Randers type}}, Canadian Mathematical Bulletin
  \textbf{57} (2014), no.~4, 765--779.

\bibitem{dehkordi2019huygens}
Hengameh~R Dehkordi and Alberto Saa, \emph{{Huygens’ envelope principle in
  Finsler spaces and analogue gravity}}, Classical and Quantum Gravity
  \textbf{36} (2019), no.~8, 085008.

\bibitem{dong2022isoparametric}
Peilong Dong and Yali Chen, \emph{{Isoparametric hypersurfaces and
  hypersurfaces with constant principal curvatures in Finsler spaces}},
  \href{https://arxiv.org/abs/2210.12937}{\texttt{arXiv.2210.12937}}, 2022.

\bibitem{dong2020isoparametric}
Peilong Dong and Qun He, \emph{{Isoparametric hypersurfaces of a class of
  Finsler manifolds induced by navigation problem in Minkowski spaces}},
  Differential Geometry and its Applications \textbf{68} (2020), 101581.

\bibitem{ge2012anisotropic}
Jianquan Ge and Hui Ma, \emph{{Anisotropic isoparametric hypersurfaces in
  Euclidean spaces}}, Annals of Global Analysis and Geometry \textbf{41}
  (2012), no.~3, 347--355.

\bibitem{gibbons2009stationary}
Gary~W Gibbons, Carlos A~R Herdeiro, Claude~M Warnick, and Marcus~C Werner,
  \emph{{Stationary metrics and optical Zermelo-Randers-Finsler geometry}},
  Physical Review D \textbf{79} (2009), no.~4, 044022.

\bibitem{gibbons2011geometry}
Gary~W Gibbons and Claude~M Warnick, \emph{{The geometry of sound rays in a
  wind}}, Contemporary Physics \textbf{52} (2011), no.~3, 197--209.

\bibitem{he2020classifications}
Qun He, Pei~Long Dong, and Song~Ting Yin, \emph{{Classifications of
  isoparametric hypersurfaces in Randers space forms}}, Acta Mathematica
  Sinica, English Series \textbf{36} (2020), no.~9, 1049--1060.

\bibitem{he2022isoparametric}
Qun He, Xin Huang, and Peilong Dong, \emph{{Isoparametric hypersurfaces in
  conic Finsler manifolds}}, Differential Geometry and its Applications
  \textbf{84} (2022), 101937.

\bibitem{he2016isoparametric}
Qun He, Songting Yin, and Yibing Shen, \emph{{Isoparametric hypersurfaces in
  Minkowski spaces}}, Differential Geometry and its Applications \textbf{47}
  (2016), 133--158.

\bibitem{he2017isoparametric}
Qun He, SongTing Yin, and YiBing Shen, \emph{{Isoparametric hypersurfaces in
  Funk spaces}}, Science China Mathematics \textbf{60} (2017), no.~12,
  2447--2464.

\bibitem{herrera2021stationary}
J{\'o}natan Herrera and Miguel~A Javaloyes, \emph{{Stationary--Complete
  Spacetimes with non-standard splittings and pre-Randers metrics}}, Journal of
  Geometry and Physics \textbf{163} (2021), 104120.

\bibitem{huang2011geodesics}
Libing Huang and Xiaohuan Mo, \emph{{On geodesics of Finsler metrics via
  navigation problem}}, Proceedings of the American Mathematical Society
  \textbf{139} (2011), no.~8, 3015--3024.

\bibitem{javaloyes2014chern}
Miguel~A Javaloyes, \emph{{Chern connection of a pseudo-Finsler metric as a
  family of affine connections}}, Publicationes Mathematicae Debrecen
  \textbf{84} (2014), no.~1--2, 29--43.

\bibitem{javaloyes2021general}
Miguel~A Javaloyes, Enrique Pend{\'a}s-Recondo, and Miguel S{\'a}nchez,
  \emph{{A general model for wildfire propagation with wind and slope}},
  \href{https://arxiv.org/abs/2110.03364}{\texttt{arXiv.2110.03364}}, 2021.

\bibitem{javaloyes2021applications}
\bysame, \emph{{Applications of cone structures to the anisotropic rheonomic
  Huygens’ principle}}, Nonlinear Analysis \textbf{209} (2021), 112337.

\bibitem{javaloyes2014definition}
Miguel~A Javaloyes and Miguel Sanchez, \emph{{On the definition and examples of
  Finsler metrics}}, Annali della Scuola Normale Superiore di Pisa. Classe di
  Scienze \textbf{13} (2014), no.~5, 813--858.

\bibitem{javaloyes2020definition}
Miguel~A Javaloyes and Miguel S{\'a}nchez, \emph{{On the definition and
  examples of cones and Finsler spacetimes}}, Revista de la Real Academia de
  Ciencias Exactas, F{\'\i}sicas y Naturales. Serie A. Matem{\'a}ticas
  \textbf{114} (2020), no.~1, 1--46.

\bibitem{javaloyes2018some}
Miguel~A Javaloyes and Henrique Vit{\'o}rio, \emph{{Some properties of Zermelo
  navigation in pseudo-Finsler metrics under an arbitrary wind}}, Houston
  Journal of Mathematics \textbf{44} (2018), no.~4, 1147--1179.

\bibitem{laura1918sopra}
Ernesto Laura, \emph{{Sopra la propagazione di onde in un mezzo indefinito}},
  Scritti Matematici Offerti ad Enrico D'Ovidio (1918), 253--278.

\bibitem{levi1937famiglie}
Tullio Levi-Civita, \emph{{Famiglie di superficie isoparametriche
  nell'ordinario spazio Euclideo}}, Atti della Accademia Nazionale dei Lincei.
  Classe di Scienze Fisiche, Matematiche e Naturali \textbf{26} (1937),
  657--664.

\bibitem{markvorsen2016finsler}
Steen Markvorsen, \emph{{A Finsler geodesic spray paradigm for wildfire spread
  modelling}}, Nonlinear Analysis: Real World Applications \textbf{28} (2016),
  208--228.

\bibitem{miyaoka2013transnormal}
Reiko Miyaoka, \emph{{Transnormal functions on a Riemannian manifold}},
  Differential Geometry and its Applications \textbf{31} (2013), no.~1,
  130--139.

\bibitem{qian2022hypersurfaces}
Yantong Qian, Qun He, and Yali Chen, \emph{{Hypersurfaces with Constant Mean
  Curvature on Finsler manifolds}},
  \href{https://arxiv.org/abs/2203.09712}{\texttt{arXiv.2203.09712}}, 2022.

\bibitem{robles2007geodesics}
Colleen Robles, \emph{{Geodesics in Randers spaces of constant curvature}},
  Transactions of the American Mathematical Society \textbf{359} (2007), no.~4,
  1633--1651.

\bibitem{schneider1997integral}
Rolf Schneider and John~A Wieacker, \emph{{Integral geometry in Minkowski
  spaces}}, Advances in Mathematics \textbf{129} (1997), no.~2, 222--260.

\bibitem{segre1924proprieta}
Beniamino Segre, \emph{{Una proprieta caratteristica di tre sistemi
  $\infty^{1}$ di superficie}}, Atti della Accademia delle Scienze di Torino.
  Classe di Scienze Fisiche, Matematiche e Naturali \textbf{59} (1924),
  666--671.

\bibitem{segre1938famiglie}
\bysame, \emph{{Famiglie di ipersuperficie isoparametriche negli spazi euclidei
  ad un qualunque numero di dimensioni}}, Atti della Accademia Nazionale dei
  Lincei. Classe di Scienze Fisiche, Matematiche e Naturali \textbf{27} (1938),
  203--207.

\bibitem{shen2016introduction}
Yi-Bing Shen and Zhongmin Shen, \emph{{Introduction to modern Finsler
  geometry}}, World Scientific Publishing Company, 2016.

\bibitem{shen1997curvature}
Zhongmin Shen, \emph{{Curvature, distance and volume in Finsler geometry}},
  Tech. Report IHES/M/97/48, Institut des Hautes {\'E}tudes Scientifiques,
  1997.

\bibitem{shen1998finsler}
\bysame, \emph{{On Finsler geometry of submanifolds}}, Mathematische Annalen
  \textbf{311} (1998), no.~3, 549--576.

\bibitem{shen1998nonlinear}
\bysame, \emph{{The non-linear Laplacian for Finsler manifolds}}, {The theory
  of Finslerian Laplacians and applications} (Peter~L Antonelli and Bradley~C
  Lackey, eds.), Springer, 1998, pp.~187--198.

\bibitem{shen2001lectures}
\bysame, \emph{{Lectures on Finsler geometry}}, World Scientific, 2001.

\bibitem{shen2003finsler}
\bysame, \emph{{Finsler Metrics with $K=0$ and $S=0$}}, Canadian Journal of
  Mathematics \textbf{55} (2003), no.~1, 112–132.

\bibitem{somigliana1918sulle}
Carlo Somigliana, \emph{{Sulle relazione fra il principio di Huygens e l'ottica
  geometrica}}, Atti della Accademia delle Scienze di Torino. Classe di Scienze
  Fisiche, Matematiche e Naturali \textbf{54} (1918--1919), 974--979.

\bibitem{souza2004bernstein}
Marcelo Souza, Joel Spruck, and Keti Tenenblat, \emph{{A Bernstein type theorem
  on a Randers space}}, Mathematische Annalen \textbf{329} (2004), no.~2,
  291--305.

\bibitem{souza2003minimal}
Marcelo Souza and Keti Tenenblat, \emph{{Minimal surfaces of rotation in
  Finsler space with a Randers metric}}, Mathematische Annalen \textbf{325}
  (2003), no.~4, 625--642.

\bibitem{thorbergsson1999survey}
Gudlaugur Thorbergsson, \emph{{A survey on isoparametric hypersurfaces and
  their generalizations}}, {Handbook of Differential Geometry} (Franki~JE
  Dillen and Leopold~CA Verstraelen, eds.), vol.~1, Elsevier, 1999,
  pp.~963--995.

\bibitem{wang1987isoparametric}
Qi-Ming Wang, \emph{{Isoparametric functions on Riemannian manifolds. I}},
  Mathematische Annalen \textbf{277} (1987), no.~4, 639--646.

\bibitem{wu2007local}
Bing~Ye Wu, \emph{{A local rigidity theorem for minimal surfaces in Minkowski
  3-space of Randers type}}, Annals of Global Analysis and Geometry \textbf{31}
  (2007), no.~4, 375--384.

\bibitem{xu2018isoparametric}
Ming Xu, \emph{{Isoparametric hypersurfaces in a Randers sphere of constant
  flag curvature}}, Annali di Matematica Pura ed Applicata (1923-) \textbf{197}
  (2018), no.~3, 703--720.

\bibitem{xu2020some}
Ming Xu, Vladimir Matveev, Ke~Yan, and Shaoxiang Zhang, \emph{{Some geometric
  correspondences for homothetic navigation}}, Publicationes Mathematicae
  Debrecen \textbf{97} (2020), no.~3--4, 449--474.

\bibitem{xu2021isoparametric}
Ming Xu, Ju~Tan, and Na~Xu, \emph{{Isoparametric hypersurfaces induced by
  navigation in Lorentz Finsler geometry}},
  \href{https://arxiv.org/abs/2105.08900}{\texttt{arXiv.2105.08900}}, 2021.

\bibitem{yajima2009finsler}
Takahiro Yajima and Hiroyuki Nagahama, \emph{{Finsler geometry of seismic ray
  path in anisotropic media}}, Proceedings of the Royal Society A:
  Mathematical, Physical and Engineering Sciences \textbf{465} (2009),
  no.~2106, 1763--1777.

\bibitem{yoshikawa2014kropina}
Ryozo Yoshikawa and Sorin~V Sabau, \emph{{Kropina metrics and Zermelo
  navigation on Riemannian manifolds}}, Geometriae Dedicata \textbf{171}
  (2014), no.~1, 119--148.

\end{thebibliography}

\end{document}